\documentclass[leqno,12pt]{amsart}
\usepackage{amssymb} 
\theoremstyle{plain} 
\newtheorem{theorem}{\indent\sc Theorem}[section] 
\newtheorem{lemma}[theorem]{\indent\sc Lemma}

\newtheorem{proposition}[theorem]{\indent\sc Proposition}

\theoremstyle{definition} 
\newtheorem{definition}[theorem]{\indent\sc Definition}

\begin{document}

\title{Holomorphy conditions of Fuji-Suzuki coupled Painlev\'e VI system \\}

\author{By\\
Yusuke Sasano}

\renewcommand{\thefootnote}{\fnsymbol{footnote}}
\footnote[0]{2000\textit{ Mathematics Subjet Classification}.
34M55; 34M45; 58F05; 32S65.}


\begin{abstract} In this note, we give some holomorphy conditions of Fuji-Suzuki coupled Painlev\'e VI system. We also give two translation operators acting  on the constant parameter $\eta$. We note a confluence process from the Fuji-Suzuki system to the Noumi-Yamada system of type $A_5^{(1)}$.
\end{abstract}
\maketitle

\section{Introduction}
In this note, we study Fuji-Suzuki coupled Painlev\'e VI system (see \cite{FS1,KNS1,KNS2}).

Define birational and symplectic transformations $r_i \ (i=0,1,\ldots,6)$ as follows:

\begin{align}\label{holo;FS}
\begin{split}
r_0:(x_0,y_0,z_0,w_0)=&\left(-((q_1-q_2)p_1-\alpha_0)p_1,\frac{1}{p_1},q_2,p_2+p_1 \right),\\
r_1:(x_1,y_1,z_1,w_1)=&\left(\frac{1}{q_1},-(q_1p_1+\alpha_1)q_1,q_2,p_2 \right),\\
r_2:(x_2,y_2,z_2,w_2)=&\left(-((q_1-t)p_1-\alpha_2)p_1,\frac{1}{p_1},q_2,p_2 \right),\\
r_3:(x_3,y_3,z_3,w_3)=&\left(-(q_1p_1+q_2p_2-(\alpha_3-\eta)))p_1,\frac{1}{p_1},q_2p_1,\frac{p_2}{p_1} \right),\\
r_4:(x_4,y_4,z_4,w_4)=&\left(q_1,p_1,-((q_2-1)p_2-\alpha_4)p_2,\frac{1}{p_2} \right),\\
r_5:(x_5,y_5,z_5,w_5)=&\left(q_1,p_1,\frac{1}{q_2},-(p_2q_2+\alpha_5)q_2 \right),\\
r_6:(x_6,y_6,z_6,w_6)=&\left(-(XY+ZW-(\eta-\alpha_1-\alpha_5))Y,\frac{1}{Y},ZY,\frac{W}{Y} \right),
\end{split}
\end{align}
where the coordinate system $(X,Y,Z,W)$ is given by
\begin{align*}
\begin{split}
&r_5 \circ r_1:(X,Y,Z,W):=\left(\frac{1}{q_1},-(q_1p_1+\alpha_1)q_1,\frac{1}{q_2},-(p_2q_2+\alpha_5)q_2 \right).
\end{split}
\end{align*}
We note that it was difficult to find the condition $r_6$. Because this condition is a patching data on the double boundary of the variables $q_1,q_2$ in 4-dimensional complex manifold $\mathcal S$ given in the paper \cite{Sasa5}, that is, $r_5 \circ r_1:(X,Y,Z,W)=\left(\frac{1}{q_1},-(q_1p_1+\alpha_1)q_1,\frac{1}{q_2},-(p_2q_2+\alpha_5)q_2 \right)$.

There exist a polynomial $H_{FS}$, such that the Hamiltonian system
\begin{equation}\label{eq:FS}
  \frac{dq_1}{dt} =\frac{\partial H_{FS}}{\partial p_1}, \ \frac{dp_1}{dt} =-\frac{\partial H_{FS}}{\partial q_1},  \ \frac{dq_2}{dt} =\frac{\partial H_{FS}}{\partial p_2}, \ \frac{dp_2}{dt} =-\frac{\partial H_{FS}}{\partial q_2}
\end{equation}
is transformed into a polynomial Hamiltonian system under the action of each $r_i \ (i=0,1,\ldots,6)$, where a polynomial Hamiltonian $H_{FS}$ is given by 
\begin{align}\label{eq:FSH}
\begin{split}
&t(t-1)H_{FS}=\\
&H_{VI}(q_1,p_1;\alpha_2,\alpha_0+\alpha_4,\alpha_3+\alpha_5-\eta,\eta \alpha_1)+H_{VI}(q_2,p_2;\alpha_0+\alpha_2,\alpha_4,\alpha_3+\alpha_1-\eta,\eta \alpha_5)\\
&+(q_1-t)(q_2-1)\{(q_1p_1+\alpha_1)p_2+p_1(q_2p_2+\alpha_5)\},
\end{split}
\end{align}
where $q_i,p_i \ (i=1,2)$ denote unknown complex variables, and $\alpha_j,\eta \ (j=0,1,\ldots,5)$ are complex constant parameters satisfying the parameter's relation:
\begin{equation*}
\alpha_0+\alpha_1+\alpha_2+\alpha_3+\alpha_4+\alpha_5=1.
\end{equation*}
It is known that the system \eqref{eq:FS},\eqref{eq:FSH} admits affine Weyl group symmetry of type $A_5^{(1)}$ as the group of its B{\"a}cklund transformations (see \cite{FS1}).

The symbol $H_{VI}(q,p;a,b,c,d)$ denotes
\begin{equation*}
H_{VI}(q,p;a,b,c,d):=q(q-1)(q-t)p^2-\{(a-1)q(q-1)+bq(q-t)+c(q-1)(q-t)\}p+dq.
\end{equation*}
We note that the holomorphy condition $r_2$ should be read that
\begin{align*}
\begin{split}
&r_2(H_{FS} - p_1)
\end{split}
\end{align*}
is a polynomial with respect to $x_2,y_2,z_2,w_2$.

This system admits several Lax pairs (see \cite{FS1,KNS1,KNS2}).

We note that the Hamiltonian system \eqref{eq:FS},\eqref{eq:FSH} is invariant under the following diagram automorphisms $s_8,s_9,s_{10}$. With {\it the notation} $(*):=(q_1,p_1,q_2,p_2,t;\alpha_0,\alpha_1,\ldots,\alpha_5,\eta)$;

\begin{align}
\begin{split}
&s_8:(*) \rightarrow \left(\frac{q_2}{t},tp_2,\frac{q_1}{t},tp_1,\frac{1}{t};\alpha_0,\alpha_5,\alpha_4,\alpha_3,\alpha_2,\alpha_1,\eta \right),\\
&s_9:(*) \rightarrow \left(\frac{1}{q_1},-(p_1q_1+p_2q_2+\eta)q_1,\frac{q_2}{q_1},p_2q_1,\frac{1}{t};\alpha_4,\alpha_3,\alpha_2,\alpha_1,\alpha_0,\alpha_5,\eta \right),\\
&s_{10}:(*) \rightarrow \left(\frac{t}{q_2},-\frac{(p_2q_2+p_1q_1+\eta)q_2}{t},\frac{q_1}{q_2},p_1q_2,t;\alpha_2,\alpha_3,\alpha_4,\alpha_5,\alpha_0,\alpha_1,\eta \right).
\end{split}
\end{align}
We also remark that these transformations $s_8,s_9,s_{10}$ satisfy the following relations:
\begin{equation*}
s_8^2=s_9^2=1, \quad s_{10}^3=1, \quad s_{10}=s_8 \circ s_9.
\end{equation*}
We remark that we can consider the transformation $s_9$ as holomorphy condition (see \eqref{th:1.4.1}, cf. \cite{Oka}).

Finally, we will give two translation operators acting  on the constant parameter $\eta$.

\begin{proposition}
Let us define the following translation operators{\rm ; \rm}
\begin{align}
\begin{split}
&T_1:=(s_2s_{10}  s_{10}  s_1)^4, \quad T_2:=s_1T_1s_1,  \quad T_3:=s_5T_1s_5.
\end{split}
\end{align}
These translation operators $T_k \ (k=1,2,3)$ act on parameters $\alpha_i,\eta$ as follows$:$
\begin{align}
\begin{split}
T_1(\alpha_0,\alpha_1,\ldots,\alpha_5,\eta)=&(\alpha_0,\alpha_1,\ldots,\alpha_5,\eta)+(0,-1,1,0,-1,1,0),\\
T_2(\alpha_0,\alpha_1,\ldots,\alpha_5,\eta)=&(\alpha_0,\alpha_1,\ldots,\alpha_5,\eta)+(-1,1,0,0,-1,1,1),\\
T_3(\alpha_0,\alpha_1,\ldots,\alpha_5,\eta)=&(\alpha_0,\alpha_1,\ldots,\alpha_5,\eta)+(1,-1,1,0,0,-1,-1).
\end{split}
\end{align}
\end{proposition}
Here, (see \cite{FS1})
\begin{align*}
\begin{split}
&s_1:(*) \rightarrow \left(q_1+\frac{\alpha_1}{p_1},p_1,q_2,p_2,t;\alpha_0+\alpha_1,-\alpha_1,\alpha_2+\alpha_1,\alpha_3,\alpha_4,\alpha_5,\eta-\alpha_1 \right),\\
&s_2:(*) \rightarrow \left(q_1,p_1-\frac{\alpha_2}{q_1-t},q_2,p_2,t;\alpha_0,\alpha_1+\alpha_2,-\alpha_2,\alpha_3+\alpha_2,\alpha_4,\alpha_5,\eta+\alpha_2 \right),\\
&s_5:(*) \rightarrow \left(q_1,p_1,q_2+\frac{\alpha_5}{p_2},p_2,t;\alpha_0+\alpha_5,\alpha_1,\alpha_2,\alpha_3,\alpha_4+\alpha_5,-\alpha_5,\eta-\alpha_5 \right).
\end{split}
\end{align*}

In particular, two transformations $T_2,T_3$ are translation operators acting  on the constant parameter $\eta:$
\begin{align}
\begin{split}
&T_2(\eta)=\eta+1, \quad T_3(\eta)=\eta-1.
\end{split}
\end{align}

Next, we review a confluence process from the system \eqref{eq:FS},\eqref{eq:FSH} to the Noumi-Yamada system of type $A_5^{(1)}$ (cf. \cite{N1,KNS1,FS1,T3}).

For the system \eqref{eq:FS},\eqref{eq:FSH}, we make a change of parameters and variables
\begin{gather}
\begin{gathered}\label{conful1}
\alpha_0=A_0, \ \alpha_1=A_1, \ldots ,\alpha_5=A_5, \quad \eta=A_1+A_5+\frac{1}{\varepsilon},
\end{gathered}\\
\begin{gathered}\label{conful2}
t=1+{\varepsilon}T, \quad q_1=1+{\varepsilon}Q_1, \quad q_2=1+{\varepsilon}Q_2, \quad p_1=\frac{P_1}{\varepsilon}, \quad p_2=\frac{P_2}{\varepsilon},
\end{gathered}\\
\begin{gathered}\label{conful3}
H_{FS}(\varepsilon)=\varepsilon \left( H_{FS}-\frac{(\alpha_1+\alpha_5) \eta}{t(t-1)(\alpha_0+\alpha_1+\alpha_2+\alpha_3+\alpha_4+\alpha_5)}   \right)
\end{gathered}
\end{gather}
from $\alpha_0,\alpha_1,\ldots,\alpha_5,\eta,t,q_1,p_1,q_2,p_2$ to $A_0,\dots ,A_5,\varepsilon,T,Q_1,P_1,Q_2,P_2$. Then the system \eqref{eq:FS},\eqref{eq:FSH} can also be written in the new variables $T,Q_1,P_1,Q_2,P_2$ and parameters $A_0,\ldots ,A_5,\\
\varepsilon$ as a Hamiltonian system
\begin{align}\label{eq:eFS}
\begin{split}
&\frac{dQ_1}{dT} =\frac{\partial H_{FS}(\varepsilon)}{\partial P_1}, \quad \frac{dP_1}{dT} =-\frac{\partial H_{FS}(\varepsilon)}{\partial Q_1},  \quad \frac{dQ_2}{dT} =\frac{\partial H_{FS}(\varepsilon)}{\partial P_2}, \quad \frac{dP_2}{dT} =-\frac{\partial H_{FS}(\varepsilon)}{\partial Q_2}.
\end{split}
\end{align}
Here, its holomorphy conditions are given by
{\footnotesize
\begin{align}\label{holo;eFS}
\begin{split}
r_0:(x_0,y_0,z_0,w_0)=&\left(-((Q_1-Q_2)P_1-A_0)P_1,\frac{1}{P_1},Q_2,P_2+P_1 \right),\\
r_1:(x_1,y_1,z_1,w_1)=&\left(\frac{1}{Q_1},-(Q_1P_1+A_1)Q_1,Q_2,P_2 \right),\\
r_2:(x_2,y_2,z_2,w_2)=&\left(-((Q_1-T)P_1-A_2)P_1,\frac{1}{P_1},Q_2,P_2 \right),\\
\tilde{r}_3:(\tilde{x}_3,\tilde{y}_3,\tilde{z}_3,\tilde{w}_3)=&\left(-\left( \left(Q_1+\frac{1}{\varepsilon} \right)P_1+\left(Q_2+\frac{1}{\varepsilon} \right)P_2-A_3+A_1+A_5+\frac{1}{\varepsilon} \right) P_1,\frac{1}{P_1},\left(Q_2+\frac{1}{\varepsilon} \right)P_1,\frac{P_2}{P_1}  \right),\\
\tilde{r}_4:(\tilde{x}_4,\tilde{y}_4,\tilde{z}_4,\tilde{w}_4)=&\left(Q_1,P_1,-(Q_2 P_2-A_4)P_2,\frac{1}{P_2} \right),\\
r_5:(x_5,y_5,z_5,w_5)=&\left(Q_1,P_1,\frac{1}{Q_2},-(P_2Q_2+A_5)Q_2 \right),\\
r_6:(x_6,y_6,z_6,w_6)=&\left(-\left(XY+ZW-\frac{1}{\varepsilon} \right)Y,\frac{1}{Y},ZY,\frac{W}{Y} \right),
\end{split}
\end{align}}
where the coordinate system $(X,Y,Z,W)$ is given by
\begin{align*}
\begin{split}
&r_5 \circ r_1:(X,Y,Z,W):=\left(\frac{1}{Q_1},-(Q_1P_1+A_1)Q_1,\frac{1}{Q_2},-(P_2Q_2+A_5)Q_2 \right).
\end{split}
\end{align*}

This new system tends to the Noumi-Yamada  system of type $A_5^{(1)}$ as $\varepsilon \rightarrow 0$, where the Noumi-Yamada system of type $A_5^{(1)}$ is explicitly given as follows: 
\begin{align}\label{eq:NYA5}
\begin{split}
&\frac{dq_1}{dt} =\frac{\partial H_{NYA5}}{\partial p_1}, \quad \frac{dp_1}{dt} =-\frac{\partial H_{NYA5}}{\partial q_1},  \quad \frac{dq_2}{dt} =\frac{\partial H_{NYA5}}{\partial p_2}, \quad \frac{dp_2}{dt} =-\frac{\partial H_{NYA5}}{\partial q_2},\\
&t H_{NYA5}=H_{V}(q_1,p_1;\alpha_0+\alpha_4,\alpha_3+\alpha_5,\alpha_1)+H_{V}(q_2,p_2;\alpha_4,\alpha_3+\alpha_1,\alpha_5)+2(q_1-t)p_1 q_2 p_2,
\end{split}
\end{align}
where $q_i,p_i \ (i=1,2)$ denote unknown complex variables, and $\alpha_j \ (j=0,1,\ldots,5)$ are complex constant parameters satisfying the parameter's relation:
\begin{equation*}
\alpha_0+\alpha_1+\alpha_2+\alpha_3+\alpha_4+\alpha_5=1.
\end{equation*}
Here, for notational convenience, we have renamed $(Q_1,P_1,Q_2,P_2,T,A_0,\dots,A_5)$ to $(q_1,p_1,\\
q_2,p_2,t,\alpha_0,\ldots,\alpha_5)$ (which are not the same as the previous $(q_1,p_1,q_2,p_2,t,\alpha_0,\ldots,\alpha_5)$).

\begin{figure}
\label{FSA5fig3}
\unitlength 0.1in
\begin{picture}( 55.1000, 29.8000)( 15.7000,-39.8000)
%
\special{pn 8}%
\special{pa 2810 1000}%
\special{pa 1800 2610}%
\special{fp}%
\special{pa 1800 2610}%
\special{pa 6910 2610}%
\special{fp}%
%
\special{pn 8}%
\special{pa 2800 1000}%
\special{pa 5810 1000}%
\special{fp}%
%
\special{pn 20}%
\special{pa 5790 1000}%
\special{pa 6900 2610}%
\special{fp}%
%
\special{pn 8}%
\special{pa 2810 1000}%
\special{pa 3440 1540}%
\special{dt 0.045}%
\special{pa 3440 1540}%
\special{pa 3070 2180}%
\special{dt 0.045}%
%
\special{pn 8}%
\special{pa 1800 2600}%
\special{pa 3070 2170}%
\special{dt 0.045}%
%
\special{pn 8}%
\special{pa 3430 1540}%
\special{pa 5320 1540}%
\special{dt 0.045}%
\special{pa 5320 1540}%
\special{pa 5790 1000}%
\special{dt 0.045}%
%
\special{pn 8}%
\special{pa 5330 1530}%
\special{pa 5680 2190}%
\special{dt 0.045}%
\special{pa 5680 2190}%
\special{pa 6890 2610}%
\special{dt 0.045}%
%
\special{pn 8}%
\special{pa 3080 2160}%
\special{pa 5650 2160}%
\special{dt 0.045}%
%
\special{pn 8}%
\special{pa 1800 2600}%
\special{pa 1570 2750}%
\special{fp}%
\special{sh 1}%
\special{pa 1570 2750}%
\special{pa 1638 2730}%
\special{pa 1616 2722}%
\special{pa 1616 2698}%
\special{pa 1570 2750}%
\special{fp}%
%
\special{pn 8}%
\special{pa 6900 2610}%
\special{pa 7080 2800}%
\special{fp}%
\special{sh 1}%
\special{pa 7080 2800}%
\special{pa 7050 2738}%
\special{pa 7044 2762}%
\special{pa 7020 2766}%
\special{pa 7080 2800}%
\special{fp}%
\put(16.7000,-29.2000){\makebox(0,0)[lb]{$Y_3$}}%
\put(67.4000,-29.0000){\makebox(0,0)[lb]{$W_4$}}%
\put(36.2000,-24.8000){\makebox(0,0)[lb]{$C_0$}}%
%
\special{pn 20}%
\special{pa 2050 2240}%
\special{pa 3230 1900}%
\special{fp}%
%
\special{pn 20}%
\special{pa 3440 1540}%
\special{pa 5320 1540}%
\special{fp}%
\put(23.4000,-20.8000){\makebox(0,0)[lb]{$C_2$}}%
\put(35.4000,-18.3000){\makebox(0,0)[lb]{$\tilde{C}_3$}}%
\put(63.2000,-17.6000){\makebox(0,0)[lb]{$\tilde{C}_4$}}%
\put(39.6000,-15.1000){\makebox(0,0)[lb]{$C_6$}}%
%
\special{pn 20}%
\special{pa 3860 2610}%
\special{pa 4400 1530}%
\special{fp}%
%
\special{pn 20}%
\special{pa 2810 1820}%
\special{pa 6000 1820}%
\special{fp}%
\put(26.7000,-33.4000){\makebox(0,0)[lb]{$\tilde{C}_3=\{(X_3,Y_3,Z_3,W_3)|X_3=Z_3=-\frac{1}{\varepsilon},Y_3=0\}$}}%
\put(30.4000,-35.9000){\makebox(0,0)[lb]{$\cup \{(X_4,Y_4,Z_4,W_4)|X_4=Z_4=-\frac{1}{\varepsilon},W_4=0\} \cong {\mathbb P}^1,$}}%
\put(26.8000,-38.8000){\makebox(0,0)[lb]{$\tilde{C}_4=\{(X_4,Y_4,Z_4,W_4)|Y_4=Z_4=W_4=0\}$}}%
\put(30.5000,-41.5000){\makebox(0,0)[lb]{$\cup \{(X_7,Y_7,Z_7,W_7)|Y_7=Z_7=W_7=0\} \cong {\mathbb P}^1.$}}%
\put(28.4000,-30.6000){\makebox(0,0)[lb]{{\bf Fuzuki-Suzuki system with patameter $\varepsilon$}}}%
\end{picture}%
\caption{This figure denotes the boundary divisor ${\mathcal H}$ of ${\mathcal S}$ (see \eqref{coordinate systems}). The bold lines $C_i \ i=0,2,6$ (see \eqref{accessible singular loci},\eqref{coordinate systems}) and $\tilde{C}_j \ j=3,4$ in ${\mathcal H}$ denote the accessible singular loci of the system \eqref{eq:eFS}.}
\end{figure}

\begin{figure}
\unitlength 0.1in
\begin{picture}( 55.1000, 42.4000)( 14.4000,-44.8000)
%
\special{pn 8}%
\special{pa 2680 920}%
\special{pa 1670 2530}%
\special{fp}%
\special{pa 1670 2530}%
\special{pa 6780 2530}%
\special{fp}%
%
\special{pn 8}%
\special{pa 2670 920}%
\special{pa 5680 920}%
\special{fp}%
%
\special{pn 20}%
\special{pa 5660 920}%
\special{pa 6770 2530}%
\special{fp}%
%
\special{pn 8}%
\special{pa 2680 920}%
\special{pa 3310 1460}%
\special{dt 0.045}%
\special{pa 3310 1460}%
\special{pa 2940 2100}%
\special{dt 0.045}%
%
\special{pn 8}%
\special{pa 1670 2520}%
\special{pa 2940 2090}%
\special{dt 0.045}%
%
\special{pn 8}%
\special{pa 3300 1460}%
\special{pa 5190 1460}%
\special{dt 0.045}%
\special{pa 5190 1460}%
\special{pa 5660 920}%
\special{dt 0.045}%
%
\special{pn 8}%
\special{pa 5200 1450}%
\special{pa 5550 2110}%
\special{dt 0.045}%
\special{pa 5550 2110}%
\special{pa 6760 2530}%
\special{dt 0.045}%
%
\special{pn 8}%
\special{pa 2950 2080}%
\special{pa 5520 2080}%
\special{dt 0.045}%
%
\special{pn 8}%
\special{pa 1670 2520}%
\special{pa 1440 2670}%
\special{fp}%
\special{sh 1}%
\special{pa 1440 2670}%
\special{pa 1508 2650}%
\special{pa 1486 2642}%
\special{pa 1486 2618}%
\special{pa 1440 2670}%
\special{fp}%
%
\special{pn 8}%
\special{pa 6770 2530}%
\special{pa 6950 2720}%
\special{fp}%
\special{sh 1}%
\special{pa 6950 2720}%
\special{pa 6920 2658}%
\special{pa 6914 2682}%
\special{pa 6890 2686}%
\special{pa 6950 2720}%
\special{fp}%
\put(15.4000,-28.4000){\makebox(0,0)[lb]{$Y_3$}}%
\put(66.1000,-28.2000){\makebox(0,0)[lb]{$W_4$}}%
\put(34.9000,-24.0000){\makebox(0,0)[lb]{$C_0$}}%
%
\special{pn 20}%
\special{pa 1920 2160}%
\special{pa 3100 1820}%
\special{fp}%
%
\special{pn 20}%
\special{pa 3310 1460}%
\special{pa 5190 1460}%
\special{fp}%
\put(22.1000,-20.0000){\makebox(0,0)[lb]{$C_2$}}%
\put(61.9000,-16.8000){\makebox(0,0)[lb]{$\tilde{C}_4$}}%
\put(38.3000,-14.3000){\makebox(0,0)[lb]{$\tilde{C}_3^{(3)}$}}%
%
\special{pn 20}%
\special{pa 3730 2530}%
\special{pa 4270 1450}%
\special{fp}%
\put(20.7000,-33.5000){\makebox(0,0)[lb]{$\tilde{C}_3^{(1)}=\{(X_6,Y_6,Z_6,W_6)|X_6=Y_6=W_6=0\}$}}%
\put(24.4000,-36.2000){\makebox(0,0)[lb]{$\cup \{(X_8,Y_8,Z_8,W_8)|X_8=Y_8=W_8=0\} \cong {\mathbb P}^1,$}}%
%
\special{pn 20}%
\special{pa 2700 930}%
\special{pa 3310 1460}%
\special{fp}%
%
\special{pn 20}%
\special{pa 5200 1470}%
\special{pa 5550 2080}%
\special{fp}%
\put(20.7000,-38.8000){\makebox(0,0)[lb]{$\tilde{C}_3^{(2)}=\{(X_{11},Y_{11},Z_{11},W_{11})|Y_{11}=Z_{11}=W_{11}=0\}$}}%
\put(24.4000,-41.5000){\makebox(0,0)[lb]{$\cup \{(X_9,Y_9,Z_9,W_9)|Y_9=Z_9=W_9=0\} \cong {\mathbb P}^1,$}}%
\put(20.7000,-43.8000){\makebox(0,0)[lb]{$\tilde{C}_3^{(3)}=\{(X_8,Y_8,Z_8,W_8)|X_8=Z_8=W_8=0\}$}}%
\put(24.4000,-46.5000){\makebox(0,0)[lb]{$\cup \{(X_9,Y_9,Z_9,W_9)|X_9=Y_9=Z_9=0\} \cong {\mathbb P}^1$}}%
\put(30.4000,-12.3000){\makebox(0,0)[lb]{$\tilde{C}_3^{(1)}$}}%
\put(54.6000,-19.3000){\makebox(0,0)[lb]{$\tilde{C}_3^{(2)}$}}%
\put(22.9000,-30.4000){\makebox(0,0)[lb]{{\bf Noumi-Yamada system of type $A_5^{(1)}$ }}}%
\put(45.0000,-6.6000){\makebox(0,0)[lb]{$\tilde{C}_3 \rightarrow C_6$}}%
\put(38.3000,-6.6000){\makebox(0,0)[lb]{$\varepsilon \rightarrow 0$}}%
\put(25.0000,-4.1000){\makebox(0,0)[lb]{{\bf Confluence process}}}%
%
\special{pn 20}%
\special{pa 4430 240}%
\special{pa 4430 840}%
\special{fp}%
\special{sh 1}%
\special{pa 4430 840}%
\special{pa 4450 774}%
\special{pa 4430 788}%
\special{pa 4410 774}%
\special{pa 4430 840}%
\special{fp}%
\end{picture}%
\label{FSA5fig4}
\caption{This figure denotes the boundary divisor ${\mathcal H}$ of ${\mathcal S}$ (see \eqref{coordinate systems}). The bold lines $C_i \ i=0,2$ (see \eqref{accessible singular loci},\eqref{coordinate systems}) and $\tilde{C}_3^{(1)},\tilde{C}_3^{(2)},\tilde{C}_3^{(3)},\tilde{C}_4$ in ${\mathcal H}$ denote the accessible singular loci of the system \eqref{eq:NYA5}.}
\end{figure}

The symbol $H_{V}(q,p;a,b,c)$ denotes
\begin{equation*}
H_{V}(q,p;a,b,c):=q(q-t)p(p+1)+a t p+b q p+c q(p+1).
\end{equation*}
Its holomorphy conditions are given by $r_j \ (j=0,1,2,5), \tilde{r}_4$ (given in \eqref{holo;eFS}) and $\tilde{r}_3^{(0)},\tilde{r}_3^{(1)},\tilde{r}_3^{(2)},\tilde{r}_3^{(3)};$ (see Figure 2)
\begin{align}
\begin{split}
&r_0:(x_0,y_0,z_0,w_0)=\left(-((q_1-q_2)p_1-\alpha_0)p_1,\frac{1}{p_1},q_2,p_2+p_1 \right),\\
&r_1:(x_1,y_1,z_1,w_1)=\left(\frac{1}{q_1},-(q_1p_1+\alpha_1)q_1,q_2,p_2 \right),\\
&r_2:(x_2,y_2,z_2,w_2)=\left(-((q_1-t)p_1-\alpha_2)p_1,\frac{1}{p_1},q_2,p_2 \right),\\
&\tilde{r}_4:(\tilde{x}_4,\tilde{y}_4,\tilde{z}_4,\tilde{w}_4)=\left(q_1,p_1,-(q_2 p_2-\alpha_4)p_2,\frac{1}{p_2} \right),\\
&r_5:(x_5,y_5,z_5,w_5)=\left(q_1,p_1,\frac{1}{q_2},-(p_2q_2+\alpha_5)q_2 \right),\\
&\tilde{r}_3^{(0)}:(\tilde{x}_3^{(0)},\tilde{y}_3^{(0)},\tilde{z}_3^{(0)},\tilde{w}_3^{(0)})=\left(\frac{1}{q_1},-((p_1+p_2+1)q_1+\alpha_3)q_1,q_2-q_1,p_2 \right),\\
&\tilde{r}_3^{(1)}:(\tilde{x}_3^{(1)},\tilde{y}_3^{(1)},\tilde{z}_3^{(1)},\tilde{w}_3^{(1)})=\left(x_1,y_1+\frac{\alpha_1-\alpha_3}{x_1}-\frac{w_1+1}{x_1^2},z_1-\frac{1}{x_1},w_1 \right),\\
&\tilde{r}_3^{(2)}:(\tilde{x}_3^{(2)},\tilde{y}_3^{(2)},\tilde{z}_3^{(2)},\tilde{w}_3^{(2)})=\left(x_5-\frac{1}{z_5},y_5,z_5,w_5+\frac{\alpha_5-\alpha_3}{z_5}-\frac{y_5+1}{z_5^2} \right),\\
&\tilde{r}_3^{(3)}:(\tilde{x}_3^{(3)},\tilde{y}_3^{(3)},\tilde{z}_3^{(3)},\tilde{w}_3^{(3)})=\left(X,Y-W+\frac{2Z W-\alpha_3+\alpha_1+\alpha_5}{X}-\frac{1}{X^2},\frac{Z-X}{X^2},W X^2 \right),
\end{split}
\end{align}
where the coordinate system $(X,Y,Z,W)$ is given by
\begin{align*}
\begin{split}
&r_5 \circ r_1:(X,Y,Z,W):=\left(\frac{1}{q_1},-(q_1p_1+\alpha_1)q_1,\frac{1}{q_2},-(p_2q_2+\alpha_5)q_2 \right).
\end{split}
\end{align*}
The Noumi-Yamada system of type $A_5^{(1)}$ can be characterized by four pairs of holomorphy conditions;
\begin{align}
\begin{split}
&\{r_0,r_1,r_2,\tilde{r}_4,r_5,\tilde{r}_3^{(0)} \}, \quad \{r_0,r_1,r_2,\tilde{r}_4,r_5,\tilde{r}_3^{(1)} \},\\
&\{r_0,r_1,r_2,\tilde{r}_4,r_5,\tilde{r}_3^{(2)} \}, \quad \{r_0,r_1,r_2,\tilde{r}_4,r_5,\tilde{r}_3^{(3)} \}.
\end{split}
\end{align}

We remark that by making a change of variables $(q_i,p_i)$ and $\alpha_j$, the following transformation $\tilde{s}_3^{(1)}$ associated with $\tilde{r}_3^{(1)}$ becomes a B{\"a}cklund transformation:
\begin{align}
\begin{split}
        &\tilde{s}_3^{(1)}: (x_1,y_1,z_1,w_1,t;\alpha_0,\alpha_1,\ldots,\alpha_5) \rightarrow \\
&\left(-x_1,-\left(y_1+\frac{\alpha_1-\alpha_3}{x_1}-\frac{w_1+1}{x_1^2} \right),z_1-\frac{1}{x_1},w_1,-t;\alpha_4,\alpha_3,\alpha_2,\alpha_1,\alpha_0,\alpha_5 \right),
        \end{split}
        \end{align}
and the following transformation $\tilde{s}_3^{(2)}$ associated with $\tilde{r}_3^{(2)}$ becomes a B{\"a}cklund transformation:
\begin{align}
\begin{split}
        &\tilde{s}_3^{(2)}: (x_5,y_5,z_5,w_5,t;\alpha_0,\alpha_1,\ldots,\alpha_5) \rightarrow \\
&\left(x_5-\frac{1}{z_5}-t,y_5,-z_5,-\left(w_5+\frac{\alpha_5-\alpha_3}{z_5}-\frac{y_5+1}{z_5^2} \right),-t;\alpha_2,\alpha_1,\alpha_0,\alpha_5,\alpha_4,\alpha_3 \right).
        \end{split}
        \end{align}
Pulling back a diagram automorphism $\pi_1;$
{\small
\begin{align}
\begin{split}
        &\pi_1: (q_1,p_1,q_2,p_2,t;\alpha_0,\alpha_1,\ldots,\alpha_5) \rightarrow  \left(-q_1,-(p_1+p_2+1),q_2-q_1,p_2,-t;\alpha_4,\alpha_3,\alpha_2,\alpha_1,\alpha_0,\alpha_5 \right)
        \end{split}
        \end{align}}
by the birational transformation $r_1$, we can obtain $\tilde{s}_3^{(1)}$, and  a diagram automorphism $\pi_2;$
{\footnotesize
\begin{align}
\begin{split}
        &\pi_2: (q_1,p_1,q_2,p_2,t;\alpha_0,\alpha_1,\ldots,\alpha_5) \rightarrow  \left(q_1-q_2-t,p_1,-q_2,-(p_2+p_1+1),-t;\alpha_2,\alpha_1,\alpha_0,\alpha_5,\alpha_4,\alpha_3 \right)
        \end{split}
        \end{align}}
by the birational transformation $r_5$, we can obtain $\tilde{s}_3^{(2)}$.

The system \eqref{eq:NYA5} has the following invariant divisors\rm{:\rm}
\begin{center}
\begin{tabular}{|c|c|c|} \hline
parameter's relation & $f_i$ \\ \hline
$\alpha_0=0$ & $f_0:=q_1-q_2$  \\ \hline
$\alpha_1=0$ & $f_1:=p_1$  \\ \hline
$\alpha_2=0$ & $f_2:=q_1-t$  \\ \hline
$\alpha_3=0$ & $f_3:=p_1+p_2+1$  \\ \hline
$\alpha_4=0$ & $f_4:=q_2$  \\ \hline
$\alpha_5=0$ & $f_5:=p_2$  \\ \hline
\end{tabular}
\end{center}
We note that when $\alpha_1=0$, we see that the system \eqref{eq:NYA5} admits a particular solution $p_1=0$, and when $\alpha_3=0$, after we make the birational and symplectic transformation:
\begin{equation*}
x_3=q_1, \ y_3=p_1+p_2+1, \ z_3=q_2-q_1, \ w_3=p_2
\end{equation*}
we see that the system \eqref{eq:NYA5} admits a particular solution $y_3=0$.

The B{\"a}cklund transformations of the system of type $A_5^{(1)}$ satisfy
\begin{equation*}
s_i(g)=g+\frac{\alpha_i}{f_i}\{f_i,g\}+\frac{1}{2!} \left(\frac{\alpha_i}{f_i} \right)^2 \{f_i,\{f_i,g\} \}+\cdots \quad (g \in {\mathbb C}(t)[q_1,p_1,q_2,p_2]),
\end{equation*}
where $\{,\}$ is the Poisson bracket such that $\{p_i,q_j\}={\delta}_{ij}$ (see \cite{N2}).

Since these B{\"a}cklund transformations have Lie theoretic origin, similarity reduction of a Drinfeld-Sokolov hierarchy admits such a B{\"a}cklund symmetry.

Finally, we list some holomorphy conditions of the system \eqref{eq:NYA5}.

{\bf Hamiltonian $H_1=r_{1}(H_{NYA5}), \ r_{1}:x=\frac{1}{q_1}, \ y=-(p_1 q_1+\alpha_1)q_1, \ z=q_2, \ w=p_2$}

\begin{align*}
&r_0^{1}:x_0=q_1, \ y_0=p_1-\frac{\alpha_0 q_2}{q_1 q_2-1}, \ z_0=q_2, \ w_0=p_2-\frac{\alpha_0 q_1}{q_1 q_2-1}, \\
&r_1^{1}:x_1=\frac{1}{q_1}, \ y_1=-(p_1 q_1+\alpha_1)q_1, \ z_1=q_2, \ w_1=p_2, \\
&r_2^{1}:x_2=-((q_1-1/t)p_1-\alpha_2)p_1, \ y_2=\frac{1}{p_1}, \ z_2=q_2 \ w_2=p_2, \\
&r_3^{1}:x_3=q_1, \ y_3=p_1-\frac{\alpha_3-\alpha_1}{q_1}-\frac{p_2+1}{q_1^2}, \ z_3=q_2-\frac{1}{q_1}, \ w_3=p_2, \\
&r_4^{1}:x_4=q_1, \ y_4=p_1, \ z_4=-(q_2 p_2-\alpha_4)p_2, \ w_4=\frac{1}{p_2}, \\
&r_5^{1}:x_5=q_1, \ y_5=p_1, \ z_5=\frac{1}{q_2}, \ w_5=-(p_2q_2+\alpha_5)q_2, 
\end{align*}
where $r_2^{1} \left(H_1+\frac{p_1}{t^2} \right)$. Here, for notational convenience, we have renamed $(x,y,z,w)$ to $(q_1,p_1,q_2,p_2)$ (which are not the same as the previous $(q_1,p_1,q_2,p_2)$).

{\bf Hamiltonian $H_5=r_{5}(H_{NYA5}), \ r_{5}:x=q_1, \ y=p_1, \ z=\frac{1}{q_2}, \ w=-(p_2q_2+\alpha_5)q_2$}

\begin{align*}
&r_0^{5}:x_0=q_1, \ y_0=p_1-\frac{\alpha_0 q_2}{q_1 q_2-1}, \ z_0=q_2, \ w_0=p_2-\frac{\alpha_0 q_1}{q_1 q_2-1}, \\
&r_1^{5}:x_1=\frac{1}{q_1}, \ y_1=-(p_1 q_1+\alpha_1)q_1, \ z_1=q_2, \ w_1=p_2, \\
&r_2^{5}:x_2=-((q_1-t)p_1-\alpha_2)p_1, \ y_2=\frac{1}{p_1}, \ z_2=q_2 \ w_2=p_2, \\
&r_3^{5}:x_3=q_1-\frac{1}{q_2}, \ y_3=p_1, \ z_3=q_2, \ w_3=p_2-\frac{\alpha_3-\alpha_5}{q_2}-\frac{p_1+1}{q_2^2}, \\
&r_4^{5}:x_4=q_1, \ y_4=p_1, \ z_4=\frac{1}{q_2}, \ w_4=-(p_2q_2+\alpha_4+\alpha_5)q_2, \\
&r_5^{5}:x_5=q_1, \ y_5=p_1, \ z_5=\frac{1}{q_2}, \ w_5=-(p_2q_2+\alpha_5)q_2, 
\end{align*}
where $r_2^{5} \left(H_5-p_1 \right)$.

\begin{figure}
\unitlength 0.1in
\begin{picture}( 35.8300, 23.4900)( 31.3000,-31.0000)
\put(43.8600,-29.5200){\makebox(0,0)[lb]{$H_{NYA5}$}}%
%
\special{pn 20}%
\special{pa 4382 2740}%
\special{pa 3414 2210}%
\special{fp}%
\special{sh 1}%
\special{pa 3414 2210}%
\special{pa 3462 2260}%
\special{pa 3460 2236}%
\special{pa 3482 2224}%
\special{pa 3414 2210}%
\special{fp}%
%
\special{pn 20}%
\special{pa 4774 2740}%
\special{pa 5724 2202}%
\special{fp}%
\special{sh 1}%
\special{pa 5724 2202}%
\special{pa 5656 2218}%
\special{pa 5678 2228}%
\special{pa 5676 2252}%
\special{pa 5724 2202}%
\special{fp}%
\put(31.5900,-21.5400){\makebox(0,0)[lb]{$H_1$}}%
\put(55.7700,-21.6900){\makebox(0,0)[lb]{$H_5$}}%
\put(42.8500,-10.1000){\makebox(0,0)[lb]{$H_{15}$}}%
%
\special{pn 20}%
\special{pa 3424 1892}%
\special{pa 4352 1080}%
\special{fp}%
\special{sh 1}%
\special{pa 4352 1080}%
\special{pa 4290 1108}%
\special{pa 4312 1114}%
\special{pa 4316 1138}%
\special{pa 4352 1080}%
\special{fp}%
%
\special{pn 20}%
\special{pa 5762 1920}%
\special{pa 4794 1080}%
\special{fp}%
\special{sh 1}%
\special{pa 4794 1080}%
\special{pa 4832 1138}%
\special{pa 4834 1114}%
\special{pa 4858 1108}%
\special{pa 4794 1080}%
\special{fp}%
\put(64.7000,-10.6000){\makebox(0,0)[lb]{$H_{153}$}}%
\put(36.6800,-26.1600){\makebox(0,0)[lb]{$r_1$}}%
\put(51.7400,-26.5100){\makebox(0,0)[lb]{$r_5$}}%
\put(35.2900,-15.5200){\makebox(0,0)[lb]{$r_5$}}%
\put(53.4200,-15.4100){\makebox(0,0)[lb]{$r_1$}}%
\put(55.0800,-9.2100){\makebox(0,0)[lb]{$r_3^{15}$}}%
%
\special{pn 20}%
\special{pa 5284 968}%
\special{pa 6212 968}%
\special{fp}%
\special{sh 1}%
\special{pa 6212 968}%
\special{pa 6146 948}%
\special{pa 6160 968}%
\special{pa 6146 988}%
\special{pa 6212 968}%
\special{fp}%
\put(35.6600,-29.1900){\makebox(0,0)[lb]{North}}%
\put(35.8600,-10.0900){\makebox(0,0)[lb]{South}}%
\put(64.4000,-21.7000){\makebox(0,0)[lb]{${H_5}'$}}%
%
\special{pn 20}%
\special{pa 6636 1120}%
\special{pa 6636 1936}%
\special{fp}%
\special{sh 1}%
\special{pa 6636 1936}%
\special{pa 6656 1870}%
\special{pa 6636 1884}%
\special{pa 6616 1870}%
\special{pa 6636 1936}%
\special{fp}%
\put(67.1300,-15.9200){\makebox(0,0)[lb]{$r_3^{153}$}}%
\put(60.4000,-26.8000){\makebox(0,0)[lb]{$\tilde{r}_5$}}%
\put(31.3000,-32.7000){\makebox(0,0)[lb]{$\tilde{r}_5:x=q_1, \ y=p_1-1, \ z=\frac{1}{q_2}, \ w=-(p_2 q_2+\alpha_5)q_2$}}%
%
\special{pn 20}%
\special{pa 6510 2230}%
\special{pa 5430 2820}%
\special{fp}%
\special{sh 1}%
\special{pa 5430 2820}%
\special{pa 5498 2806}%
\special{pa 5478 2794}%
\special{pa 5480 2770}%
\special{pa 5430 2820}%
\special{fp}%
\end{picture}%
\label{fig:•Ïg‹Èü'ÌŠGFSNYA5G}
\caption{Relation between Hamiltonians $H_{NYA5}$ and $H_1,H_5,H_{15},H_{153}$}
\end{figure}
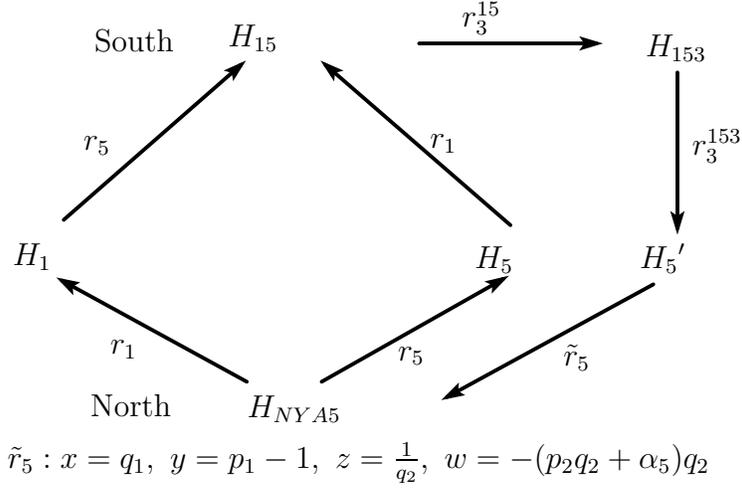

{\bf Hamiltonian $H_{15}=r_{15}(H_{NYA5}), \ r_{15}:x=\frac{1}{q_1}, \ y=-(p_1 q_1+\alpha_1)q_1, \ z=\frac{1}{q_2}, \ w=-(p_2q_2+\alpha_5)q_2$}

\begin{align*}
&r_0^{15}:x_0=-((q_1-q_2)p_1-\alpha_0)p_1, \ y_0=\frac{1}{p_1}, \ z_0=q_2, \ w_0=p_2+p_1, \\
&r_1^{15}:x_1=\frac{1}{q_1}, \ y_1=-(p_1 q_1+\alpha_1)q_1, \ z_1=q_2, \ w_1=p_2, \\
&r_2^{15}:x_2=-((q_1-1/t)p_1-\alpha_2)p_1, \ y_2=\frac{1}{p_1}, \ z_2=q_2 \ w_2=p_2, \\
&r_3^{15}:x_3=q_1, \ y_3=p_1-p_2+\frac{2q_2 p_2-(\alpha_3-\alpha_1-\alpha_5)}{q_1}-\frac{1}{q_1^2}, \ z_3=\frac{q_2-q_1}{q_1^2}, \ w_3=p_2 q_1^2, \\
&r_4^{15}:x_4=q_1, \ y_4=p_1, \ z_4=\frac{1}{q_2}, \ w_4=-(p_2q_2+\alpha_4+\alpha_5)q_2, \\
&r_5^{15}:x_5=q_1, \ y_5=p_1, \ z_5=\frac{1}{q_2}, \ w_5=-(p_2q_2+\alpha_5)q_2, 
\end{align*}
where $r_2^{15} \left(H_{15}+\frac{p_1}{t^2} \right)$.

{\bf Hamiltonian $H_{153}=r_3^{15}(H_{15})$}

{\footnotesize
\begin{align*}
&r_0^{153}:x_0=q_1, \ y_0=p_1, \ z_0=-(q_2 p_2-\alpha_0)p_2, \ w_0=\frac{1}{p_2}, \\
&r_1^{153}:x_1=q_1, \ y_1=p_1-\frac{2q_2 p_2-(\alpha_3-\alpha_1-\alpha_5)}{q_1}+\frac{1}{q_1^2}, \ z_1=q_2 q_1^2, \ w_1=\frac{p_2}{q_1^2}, \\
&r_2^{153}:x_2=-((q_1-1/t)p_1-\alpha_2)p_1, \ y_2=\frac{1}{p_1}, \ z_2=q_2 \ w_2=p_2, \\
&r_3^{153}:x_3=\frac{1}{q_1}, \ y_3=-\left( \left(p_1+\frac{p_2}{q_1^2} \right)-2(q_2 q_1+1)\frac{p_2}{q_1}+\alpha_3-\alpha_5 \right)q_1, \ z_3=(q_2 q_1+1)q_1, \ w_3=\frac{p_2}{q_1^2}, \\
&r_4^{153}:x_4=q_1, \ y_4=p_1, \ z_4=\frac{1}{q_2}, \ w_4=-(p_2q_2+\alpha_4+\alpha_5)q_2, \\
&r_5^{153}:x_5=q_1, \ y_5=p_1, \ z_5=\frac{1}{q_2}, \ w_5=-(p_2q_2+\alpha_5)q_2, 
\end{align*}}
where $r_2^{153} \left(H_{153}+\frac{p_1}{t^2} \right)$.

After we review the notion of accessible singularity in next section, we make its holomorphy conditions by resolving the accessible singularities.

\section{Accessible singularity and local index}
Let us review the notion of {\it accessible singularity}. Let $B$ be a connected open domain in $\mathbb C$ and $\pi : {\mathcal W} \longrightarrow B$ a smooth proper holomorphic map. We assume that ${\mathcal H} \subset {\mathcal W}$ is a normal crossing divisor which is flat over $B$. Let us consider a rational vector field $\tilde v$ on $\mathcal W$ satisfying the condition
\begin{equation*}
\tilde v \in H^0({\mathcal W},\Theta_{\mathcal W}(-\log{\mathcal H})({\mathcal H})).
\end{equation*}
Fixing $t_0 \in B$ and $P \in {\mathcal W}_{t_0}$, we can take a local coordinate system $(x_1,\ldots ,x_n)$ of ${\mathcal W}_{t_0}$ centered at $P$ such that ${\mathcal H}_{\rm smooth \rm}$ can be defined by the local equation $x_1=0$.
Since $\tilde v \in H^0({\mathcal W},\Theta_{\mathcal W}(-\log{\mathcal H})({\mathcal H}))$, we can write down the vector field $\tilde v$ near $P=(0,\ldots ,0,t_0)$ as follows:
\begin{equation*}
\tilde v= \frac{\partial}{\partial t}+g_1 
\frac{\partial}{\partial x_1}+\frac{g_2}{x_1} 
\frac{\partial}{\partial x_2}+\cdots+\frac{g_n}{x_1} 
\frac{\partial}{\partial x_n}.
\end{equation*}
This vector field defines the following system of differential equations
\begin{equation}\label{39}
\frac{dx_1}{dt}=g_1(x_1,\ldots,x_n,t),\ \frac{dx_2}{dt}=\frac{g_2(x_1,\ldots,x_n,t)}{x_1},\cdots, \frac{dx_n}{dt}=\frac{g_n(x_1,\ldots,x_n,t)}{x_1}.
\end{equation}
Here $g_i(x_1,\ldots,x_n,t), \ i=1,2,\ldots ,n,$ are holomorphic functions defined near $P$.

\begin{definition}\label{Def1}
With the above notation, assume that the rational vector field $\tilde v$ on $\mathcal W$ satisfies the condition
$$
(A) \quad \tilde v \in H^0({\mathcal W},\Theta_{\mathcal W}(-\log{\mathcal H})({\mathcal H})).
$$
We say that $\tilde v$ has an {\it accessible singularity} at $P=(0,\dots ,0,t_0)$ if
\begin{equation}
\boxed{%
x_1=0 \ {\rm and \rm} \ g_i(0,\ldots,0,t_0)=0 \ {\rm for \rm} \ {\rm every \rm} \ i, \ 2 \leq i \leq n.
}%
\end{equation}
\end{definition}

If $P \in {\mathcal H}_{{\rm smooth \rm}}$ is not an accessible singularity, all solutions of the ordinary differential equation passing through $P$ are vertical solutions, that is, the solutions are contained in the fiber ${\mathcal W}_{t_0}$ over $t=t_0$. If $P \in {\mathcal H}_{\rm smooth \rm}$ is an accessible singularity, there may be a solution of \eqref{39} which passes through $P$ and goes into the interior ${\mathcal W}-{\mathcal H}$ of ${\mathcal W}$.

\section{Construction of the holomorphy conditions}
In this section, we will give the holomorphy conditions $r_i \ (i=0,1,\ldots,6)$ by resolving some accessible singular loci of the system \eqref{eq:FS},\eqref{eq:FSH}.

In order to consider the singularity analysis for the system \eqref{eq:FS},\eqref{eq:FSH} as a compactification of ${\mathbb C}^4$ which is the phase space of the system \eqref{eq:FS},\eqref{eq:FSH}, we take 4-dimensional complex manifold $\mathcal S$ given in the paper \cite{Sasa5}. This manifold can be considered as a generalization of the Hirzebruch surface.

We easily see that the rational vector field $\tilde v$ associated with the system \eqref{eq:FS},\eqref{eq:FSH} satisfies the condition:
$$
\tilde v \in H^0({\mathcal S},\Theta_{\mathcal S}(-\log{\mathcal H})({\mathcal H})).
$$

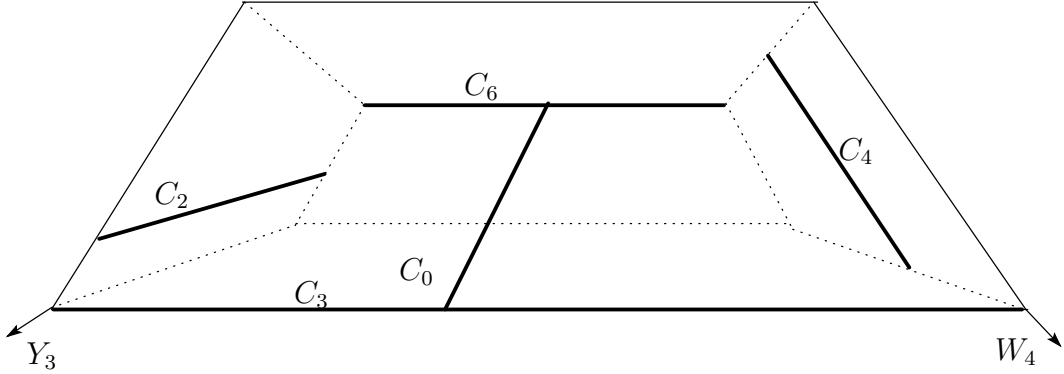
\begin{figure}
\unitlength 0.1in
\begin{picture}( 55.1000, 18.0000)( 15.7000,-28.0000)
%
\special{pn 8}%
\special{pa 2810 1000}%
\special{pa 1800 2610}%
\special{fp}%
\special{pa 1800 2610}%
\special{pa 6910 2610}%
\special{fp}%
%
\special{pn 8}%
\special{pa 2800 1000}%
\special{pa 5810 1000}%
\special{fp}%
%
\special{pn 8}%
\special{pa 5790 1000}%
\special{pa 6900 2610}%
\special{fp}%
%
\special{pn 8}%
\special{pa 2810 1000}%
\special{pa 3440 1540}%
\special{dt 0.045}%
\special{pa 3440 1540}%
\special{pa 3070 2180}%
\special{dt 0.045}%
%
\special{pn 8}%
\special{pa 1800 2600}%
\special{pa 3070 2170}%
\special{dt 0.045}%
%
\special{pn 8}%
\special{pa 3430 1540}%
\special{pa 5320 1540}%
\special{dt 0.045}%
\special{pa 5320 1540}%
\special{pa 5790 1000}%
\special{dt 0.045}%
%
\special{pn 8}%
\special{pa 5330 1530}%
\special{pa 5680 2190}%
\special{dt 0.045}%
\special{pa 5680 2190}%
\special{pa 6890 2610}%
\special{dt 0.045}%
%
\special{pn 8}%
\special{pa 3080 2160}%
\special{pa 5650 2160}%
\special{dt 0.045}%
%
\special{pn 8}%
\special{pa 1800 2600}%
\special{pa 1570 2750}%
\special{fp}%
\special{sh 1}%
\special{pa 1570 2750}%
\special{pa 1638 2730}%
\special{pa 1616 2722}%
\special{pa 1616 2698}%
\special{pa 1570 2750}%
\special{fp}%
%
\special{pn 8}%
\special{pa 6900 2610}%
\special{pa 7080 2800}%
\special{fp}%
\special{sh 1}%
\special{pa 7080 2800}%
\special{pa 7050 2738}%
\special{pa 7044 2762}%
\special{pa 7020 2766}%
\special{pa 7080 2800}%
\special{fp}%
\put(16.7000,-29.2000){\makebox(0,0)[lb]{$Y_3$}}%
\put(67.4000,-29.0000){\makebox(0,0)[lb]{$W_4$}}%
%
\special{pn 20}%
\special{pa 1810 2610}%
\special{pa 6880 2610}%
\special{fp}%
\put(36.2000,-24.8000){\makebox(0,0)[lb]{$C_0$}}%
%
\special{pn 20}%
\special{pa 2050 2240}%
\special{pa 3230 1900}%
\special{fp}%
%
\special{pn 20}%
\special{pa 6290 2390}%
\special{pa 5550 1280}%
\special{fp}%
%
\special{pn 20}%
\special{pa 3440 1540}%
\special{pa 5320 1540}%
\special{fp}%
\put(23.4000,-20.8000){\makebox(0,0)[lb]{$C_2$}}%
\put(30.7000,-26.0000){\makebox(0,0)[lb]{$C_3$}}%
\put(59.2000,-18.5000){\makebox(0,0)[lb]{$C_4$}}%
\put(39.6000,-15.1000){\makebox(0,0)[lb]{$C_6$}}%
%
\special{pn 20}%
\special{pa 3860 2610}%
\special{pa 4400 1530}%
\special{fp}%
\end{picture}%
\label{FSA5fig1}
\caption{This figure denotes the boundary divisor ${\mathcal H}$ of ${\mathcal S}$. The bold lines $C_i \ i=0,2,3,4,6$ in ${\mathcal H}$ denote the accessible singular loci of the system \eqref{eq:FS},\eqref{eq:FSH}.}
\end{figure}

\begin{lemma}\label{lem1}
The rational vector field $\tilde v$ associated with the system \eqref{eq:FS},\eqref{eq:FSH}  has the following accessible singular loci $C_i \cong {\mathbb P}^1 \ (i=0,2,3,4,6)$ \rm{(see Figure 4)}$:$
\begin{equation}\label{accessible singular loci}
  \left\{
  \begin{aligned}
    C_0=&\{(X_3,Y_3,Z_3,W_3)|X_3=Z_3,Y_3=0,W_3=-1\}\\
    &\cup \{(X_8,Y_8,Z_8,W_8)|X_8=Z_8,Y_8=0,W_8=-1\} \cong {\mathbb P}^1,\\
    C_2=&\{(X_3,Y_3,Z_3,W_3)|X_3=t,Y_3=0,W_3=0\},\\
    &\cup \{(X_{10},Y_{10},Z_{10},W_{10})|X_{10}=t,Y_{10}=0,W_{10}=0\} \cong {\mathbb P}^1,\\
    C_3=&\{(X_3,Y_3,Z_3,W_3)|X_3=Y_3=Z_3=0\}\\
&\cup \{(X_4,Y_4,Z_4,W_4)|X_4=Z_4=W_4=0\} \cong {\mathbb P}^1,\\
    C_4=&\{(X_4,Y_4,Z_4,W_4)|Y_4=0,Z_4=1,W_4=0\},\\
    &\cup \{(X_7,Y_7,Z_7,W_7)|Y_7=0,Z_7=1,W_7=0\} \cong {\mathbb P}^1,\\
    C_6=&\{(X_8,Y_8,Z_8,W_8)|X_8=Y_8=Z_8=0\}\\
&\cup \{(X_9,Y_9,Z_9,W_9)|X_9=Z_9=W_9=0\} \cong {\mathbb P}^1.
   \end{aligned}
  \right. 
\end{equation}
\end{lemma}

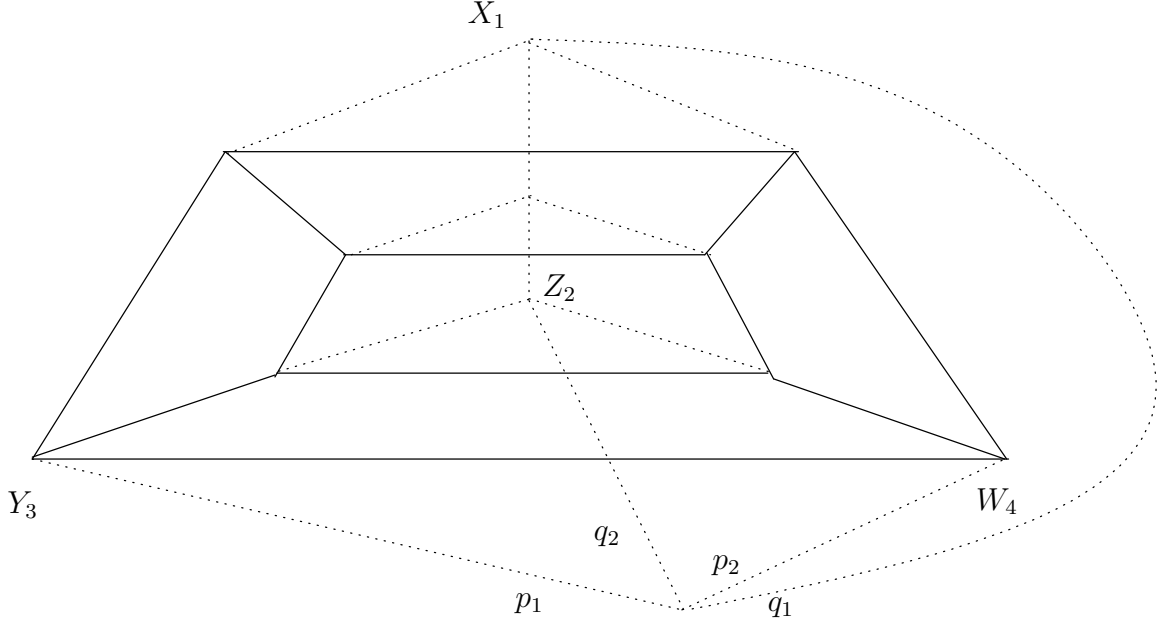
\begin{figure}
\unitlength 0.1in
\begin{picture}( 60.1200, 32.2000)( 16.7000,-34.0000)
\put(16.7000,-29.2000){\makebox(0,0)[lb]{$Y_3$}}%
%
\special{pn 8}%
\special{pa 2810 1000}%
\special{pa 1800 2610}%
\special{fp}%
\special{pa 1800 2610}%
\special{pa 6910 2610}%
\special{fp}%
%
\special{pn 8}%
\special{pa 2800 1000}%
\special{pa 5810 1000}%
\special{fp}%
%
\special{pn 8}%
\special{pa 5790 1000}%
\special{pa 6900 2610}%
\special{fp}%
%
\special{pn 8}%
\special{pa 2810 1000}%
\special{pa 3440 1540}%
\special{fp}%
\special{pa 3440 1540}%
\special{pa 3070 2180}%
\special{fp}%
%
\special{pn 8}%
\special{pa 1800 2600}%
\special{pa 3070 2170}%
\special{fp}%
%
\special{pn 8}%
\special{pa 3430 1540}%
\special{pa 5320 1540}%
\special{fp}%
\special{pa 5320 1540}%
\special{pa 5790 1000}%
\special{fp}%
%
\special{pn 8}%
\special{pa 5330 1530}%
\special{pa 5680 2190}%
\special{fp}%
\special{pa 5680 2190}%
\special{pa 6890 2610}%
\special{fp}%
%
\special{pn 8}%
\special{pa 3080 2160}%
\special{pa 5650 2160}%
\special{fp}%
\put(67.4000,-29.0000){\makebox(0,0)[lb]{$W_4$}}%
%
\special{pn 8}%
\special{pa 2820 1010}%
\special{pa 4420 410}%
\special{dt 0.045}%
%
\special{pn 8}%
\special{pa 4400 430}%
\special{pa 5790 990}%
\special{dt 0.045}%
%
\special{pn 8}%
\special{pa 4400 420}%
\special{pa 4400 1760}%
\special{dt 0.045}%
%
\special{pn 8}%
\special{pa 3090 2150}%
\special{pa 4410 1770}%
\special{dt 0.045}%
%
\special{pn 8}%
\special{pa 4400 1770}%
\special{pa 5680 2160}%
\special{dt 0.045}%
%
\special{pn 8}%
\special{pa 3470 1540}%
\special{pa 4410 1230}%
\special{dt 0.045}%
%
\special{pn 8}%
\special{pa 4400 1240}%
\special{pa 5350 1540}%
\special{dt 0.045}%
%
\special{pn 8}%
\special{pa 1800 2610}%
\special{pa 5200 3400}%
\special{dt 0.045}%
\special{pa 5200 3400}%
\special{pa 6890 2610}%
\special{dt 0.045}%
%
\special{pn 8}%
\special{pa 4400 1780}%
\special{pa 5210 3400}%
\special{dt 0.045}%
%
\special{pn 8}%
\special{pa 4410 410}%
\special{pa 4444 412}%
\special{pa 4476 414}%
\special{pa 4510 414}%
\special{pa 4542 416}%
\special{pa 4576 416}%
\special{pa 4608 418}%
\special{pa 4642 420}%
\special{pa 4674 420}%
\special{pa 4708 422}%
\special{pa 4740 424}%
\special{pa 4774 426}%
\special{pa 4806 428}%
\special{pa 4840 428}%
\special{pa 4872 430}%
\special{pa 4906 432}%
\special{pa 4938 434}%
\special{pa 4970 436}%
\special{pa 5004 438}%
\special{pa 5036 442}%
\special{pa 5070 444}%
\special{pa 5102 446}%
\special{pa 5134 448}%
\special{pa 5166 452}%
\special{pa 5200 454}%
\special{pa 5232 458}%
\special{pa 5264 460}%
\special{pa 5296 464}%
\special{pa 5330 468}%
\special{pa 5362 472}%
\special{pa 5394 476}%
\special{pa 5426 480}%
\special{pa 5458 484}%
\special{pa 5490 488}%
\special{pa 5522 494}%
\special{pa 5554 498}%
\special{pa 5586 504}%
\special{pa 5618 508}%
\special{pa 5650 514}%
\special{pa 5682 520}%
\special{pa 5712 526}%
\special{pa 5744 532}%
\special{pa 5776 538}%
\special{pa 5808 546}%
\special{pa 5838 552}%
\special{pa 5870 560}%
\special{pa 5900 568}%
\special{pa 5932 576}%
\special{pa 5962 584}%
\special{pa 5994 592}%
\special{pa 6024 600}%
\special{pa 6054 610}%
\special{pa 6086 620}%
\special{pa 6116 628}%
\special{pa 6146 638}%
\special{pa 6176 650}%
\special{pa 6206 660}%
\special{pa 6236 670}%
\special{pa 6266 682}%
\special{pa 6296 694}%
\special{pa 6326 706}%
\special{pa 6356 718}%
\special{pa 6386 732}%
\special{pa 6416 744}%
\special{pa 6444 758}%
\special{pa 6474 772}%
\special{pa 6502 786}%
\special{pa 6532 802}%
\special{pa 6560 816}%
\special{pa 6590 832}%
\special{pa 6618 848}%
\special{pa 6646 864}%
\special{pa 6674 882}%
\special{pa 6702 900}%
\special{pa 6730 916}%
\special{pa 6758 936}%
\special{pa 6786 954}%
\special{pa 6814 974}%
\special{pa 6842 992}%
\special{pa 6868 1012}%
\special{pa 6896 1034}%
\special{pa 6922 1054}%
\special{pa 6950 1076}%
\special{pa 6976 1098}%
\special{pa 7002 1120}%
\special{pa 7030 1144}%
\special{pa 7056 1168}%
\special{pa 7082 1192}%
\special{pa 7108 1216}%
\special{pa 7132 1242}%
\special{pa 7158 1266}%
\special{pa 7182 1292}%
\special{pa 7208 1318}%
\special{pa 7232 1344}%
\special{pa 7256 1372}%
\special{pa 7278 1398}%
\special{pa 7302 1426}%
\special{pa 7324 1452}%
\special{pa 7346 1480}%
\special{pa 7368 1508}%
\special{pa 7390 1536}%
\special{pa 7410 1564}%
\special{pa 7430 1592}%
\special{pa 7450 1620}%
\special{pa 7468 1648}%
\special{pa 7486 1678}%
\special{pa 7504 1706}%
\special{pa 7520 1734}%
\special{pa 7536 1762}%
\special{pa 7552 1792}%
\special{pa 7566 1820}%
\special{pa 7580 1848}%
\special{pa 7594 1876}%
\special{pa 7606 1904}%
\special{pa 7618 1932}%
\special{pa 7628 1960}%
\special{pa 7638 1988}%
\special{pa 7648 2014}%
\special{pa 7654 2042}%
\special{pa 7662 2068}%
\special{pa 7668 2094}%
\special{pa 7672 2120}%
\special{pa 7676 2146}%
\special{pa 7680 2172}%
\special{pa 7682 2198}%
\special{pa 7682 2222}%
\special{pa 7682 2246}%
\special{pa 7680 2270}%
\special{pa 7678 2294}%
\special{pa 7676 2318}%
\special{pa 7670 2340}%
\special{pa 7666 2362}%
\special{pa 7660 2386}%
\special{pa 7652 2408}%
\special{pa 7644 2428}%
\special{pa 7634 2450}%
\special{pa 7624 2470}%
\special{pa 7614 2492}%
\special{pa 7602 2512}%
\special{pa 7588 2532}%
\special{pa 7574 2550}%
\special{pa 7560 2570}%
\special{pa 7544 2590}%
\special{pa 7528 2608}%
\special{pa 7512 2626}%
\special{pa 7494 2644}%
\special{pa 7474 2662}%
\special{pa 7454 2680}%
\special{pa 7434 2696}%
\special{pa 7414 2714}%
\special{pa 7392 2730}%
\special{pa 7370 2746}%
\special{pa 7346 2764}%
\special{pa 7322 2778}%
\special{pa 7298 2794}%
\special{pa 7272 2810}%
\special{pa 7246 2826}%
\special{pa 7220 2840}%
\special{pa 7192 2854}%
\special{pa 7164 2870}%
\special{pa 7134 2884}%
\special{pa 7106 2898}%
\special{pa 7076 2912}%
\special{pa 7046 2926}%
\special{pa 7014 2938}%
\special{pa 6982 2952}%
\special{pa 6950 2964}%
\special{pa 6918 2978}%
\special{pa 6884 2990}%
\special{pa 6850 3002}%
\special{pa 6816 3014}%
\special{pa 6782 3026}%
\special{pa 6746 3038}%
\special{pa 6712 3050}%
\special{pa 6674 3062}%
\special{pa 6638 3074}%
\special{pa 6602 3086}%
\special{pa 6564 3096}%
\special{pa 6526 3108}%
\special{pa 6488 3118}%
\special{pa 6450 3128}%
\special{pa 6410 3140}%
\special{pa 6372 3150}%
\special{pa 6332 3160}%
\special{pa 6292 3170}%
\special{pa 6252 3180}%
\special{pa 6210 3190}%
\special{pa 6170 3200}%
\special{pa 6128 3210}%
\special{pa 6088 3220}%
\special{pa 6046 3230}%
\special{pa 6004 3240}%
\special{pa 5962 3248}%
\special{pa 5920 3258}%
\special{pa 5878 3268}%
\special{pa 5834 3276}%
\special{pa 5792 3286}%
\special{pa 5748 3294}%
\special{pa 5706 3304}%
\special{pa 5662 3314}%
\special{pa 5618 3322}%
\special{pa 5576 3332}%
\special{pa 5532 3340}%
\special{pa 5488 3348}%
\special{pa 5444 3358}%
\special{pa 5400 3366}%
\special{pa 5356 3376}%
\special{pa 5312 3384}%
\special{pa 5268 3392}%
\special{pa 5230 3400}%
\special{sp -0.045}%
\put(43.3000,-34.2000){\makebox(0,0)[lb]{$p_1$}}%
\put(47.4000,-30.7000){\makebox(0,0)[lb]{$q_2$}}%
\put(53.6000,-32.2000){\makebox(0,0)[lb]{$p_2$}}%
\put(56.5000,-34.4000){\makebox(0,0)[lb]{$q_1$}}%
\put(40.8000,-3.5000){\makebox(0,0)[lb]{$X_1$}}%
\put(44.7000,-17.8000){\makebox(0,0)[lb]{$Z_2$}}%
\end{picture}%
\label{4dimHirzebruchmanifold}
\caption{This figure denotes 4-dimensional complex manifold $\mathcal S$ (see  \cite{Sasa5}) and its boundary divisor ${\mathcal H}$.  ${\mathcal H}$ is drawn by solid line.}
\end{figure}

Here, the coordinate systems $(X_i,Y_i,Z_i,W_i) \ (i=0,1,\cdots,11)$ (see  Figure 4, cf. \cite{Sasa5}) are explicitly given by

{\small
\begin{align}\label{coordinate systems}
\begin{split}
&(X_0,Y_0,Z_0,W_0)=\left(q_1,p_1,q_2,p_2 \right),\\
&(X_1,Y_1,Z_1,W_1)=\left(\frac{1}{q_1},-(p_1q_1+\alpha_1)q_1,q_2,p_2 \right),\\
&(X_2,Y_2,Z_2,W_2)=\left(q_1,p_1,\frac{1}{q_2},-(p_2q_2+\alpha_5)q_2 \right),\\
&(X_3,Y_3,Z_3,W_3)=\left(q_1,\frac{1}{p_1},q_2,\frac{p_2}{p_1} \right),\\
&(X_4,Y_4,Z_4,W_4)=\left(q_1,\frac{p_1}{p_2},q_2,\frac{1}{p_2} \right),\\
&(X_5,Y_5,Z_5,W_5)=\left(\frac{1}{q_1},-(p_1q_1+\alpha_1)q_1,\frac{1}{q_2},-(p_2q_2+\alpha_5)q_2 \right),\\
&(X_6,Y_6,Z_6,W_6)=\left(\frac{1}{q_1},-\frac{1}{(q_1p_1+\alpha_1)q_1},q_2,-\frac{p_2}{(q_1p_1+\alpha_1)q_1} \right),\\
&(X_7,Y_7,Z_7,W_7)=\left(\frac{1}{q_1},-\frac{(p_1q_1+\alpha_1)q_1}{p_2},q_2,\frac{1}{p_2} \right),\\
&(X_8,Y_8,Z_8,W_8)=\left(\frac{1}{q_1},-\frac{1}{(p_1q_1+\alpha_1)q_1},\frac{1}{q_2},\frac{(p_2q_2+\alpha_5)q_2}{(p_1q_1+\alpha_1)q_1} \right),\\
&(X_9,Y_9,Z_9,W_9)=\left(\frac{1}{q_1},\frac{(p_1q_1+\alpha_1)q_1}{(p_2q_2+\alpha_5)q_2},\frac{1}{q_2},-\frac{1}{(p_2q_2+\alpha_5)q_2} \right),\\
&(X_{10},Y_{10},Z_{10},W_{10})=\left(q_1,\frac{1}{p_1},\frac{1}{q_2},-\frac{(p_2q_2+\alpha_5)q_2}{p_1} \right),\\
&(X_{11},Y_{11},Z_{11},W_{11})=\left(q_1,-\frac{p_1}{(p_2q_2+\alpha_5)q_2},\frac{1}{q_2},-\frac{1}{(p_2q_2+\alpha_5)q_2} \right).
\end{split}
\end{align}}

\begin{proposition}\label{prop3}
If we resolve the accessible singular loci given in Lemma \ref{lem1} by blowing-ups, then we can obtain the canonical coordinates $r_i \ (i=0,2,3,4,6)$.
\end{proposition}

{\it Proof.} By the following steps, we can resolve the accessible singular locus $C_3$.

{\bf Step 0:} Around the point $P:=\{(X_3,Y_3,Z_3,W_3)|X_3=Y_3=Z_3=W_3=0\}$, we rewrite the system \eqref{eq:FS} as follows:
\begin{align*}
\frac{d}{dt}\begin{pmatrix}
             X_3 \\
             Y_3 \\
             Z_3 \\
             W_3
             \end{pmatrix}&=\frac{1}{Y_1}\left\{\begin{pmatrix}
             \frac{2}{t-1} & -\frac{\alpha_3-\eta}{t-1} & 0 & 0  \\
             0 & \frac{1}{t-1} & 0 & 0 \\
             \frac{1}{t-1} & 0 & \frac{1}{t-1} & 0 \\
             0 & 0 & 0 & 0
             \end{pmatrix}\begin{pmatrix}
             X_3 \\
             Y_3 \\
             Z_3 \\
             W_3
             \end{pmatrix}+\cdots\right\}.
             \end{align*}

{\bf Step 1:} We blow up along the curve $C_3$.
$$
X_3^{(1)}=\frac{X_3}{Y_3}, \quad Y_3^{(1)}=Y_3, \quad Z_3^{(1)}=\frac{Z_3}{Y_3}, \quad W_3^{(1)}=W_3.
$$

{\bf Step 2:} We blow up along the surface $\{(X_3^{(1)},Y_3^{(1)},Z_3^{(1)},W_3^{(1)})|X_3^{(1)}=-Z_3^{(1)} W_3^{(1)}+(\alpha_3-\eta)\}$
$$
X_3^{(2)}=\frac{X_3^{(1)}+Z_3^{(1)} W_3^{(1)}-(\alpha_3-\eta)}{Y_3^{(1)}}, \quad Y_3^{(2)}=Y_3^{(1)}, \quad Z_3^{(2)}=Z_3^{(1)}, \quad W_3^{(2)}=W_3^{(1)}.
$$
Now, we have resolved the accessible singularity $C_3$.

By choosing a new coordinate system as
$$
(x_3,y_3,z_3,w_3)=(-X_3^{(2)},Y_3^{(2)},Z_3^{(2)},W_3^{(2)}),
$$
we can obtain the coordinate $r_3$.

Next, by the following steps, we can resolve the accessible singular locus $C_6$.

{\bf Step 0:} Around the point $Q:=\{(X_8,Y_8,Z_8,W_8)|X_8=Y_8=Z_8=W_8=0\}$, we rewrite the system \eqref{eq:FS} as follows:
\begin{align*}
\frac{d}{dt}\begin{pmatrix}
             X_8 \\
             Y_8 \\
             Z_8 \\
             W_8
             \end{pmatrix}&=\frac{1}{Y_1}\left\{\begin{pmatrix}
             \frac{2}{t(t-1)} & -\frac{\eta-(\alpha_1+\alpha_5)}{t(t-1)} & 0 & 0  \\
             0 & \frac{1}{t(t-1)} & 0 & 0 \\
             \frac{1}{t(t-1)} & 0 & \frac{1}{t(t-1)} & 0 \\
             0 & 0 & 0 & 0
             \end{pmatrix}\begin{pmatrix}
             X_8 \\
             Y_8 \\
             Z_8 \\
             W_8
             \end{pmatrix}+\cdots\right\}.
             \end{align*}

{\bf Step 1:} We blow up along the curve $C_6$.
$$
X_8^{(1)}=\frac{X_8}{Y_8}, \quad Y_8^{(1)}=Y_8, \quad Z_8^{(1)}=\frac{Z_8}{Y_8}, \quad W_8^{(1)}=W_8.
$$

{\bf Step 2:} We blow up along the surface $\{(X_8^{(1)},Y_8^{(1)},Z_8^{(1)},W_8^{(1)})|X_8^{(1)}=-Z_8^{(1)} W_8^{(1)}+\eta-(\alpha_1+\alpha_5)\}$
$$
X_8^{(2)}=\frac{X_8^{(1)}+Z_8^{(1)} W_8^{(1)}-(\eta-\alpha_1-\alpha_5)}{Y_8^{(1)}}, \quad Y_8^{(2)}=Y_8^{(1)}, \quad Z_8^{(2)}=Z_8^{(1)}, \quad W_8^{(2)}=W_8^{(1)}.
$$
Now, we have resolved the accessible singularity $C_6$. 

By choosing a new coordinate system as
$$
(x_6,y_6,z_6,w_6)=(-X_8^{(2)},Y_8^{(2)},Z_8^{(2)},W_8^{(2)}),
$$
we can obtain the coordinate $r_6$.

For the remaining accessible singular loci, the proof is similar. Collecting all the cases, we have obtained the canonical coordinates $r_i \ (i=0,2,3,4,6)$, which proves Proposition \ref{prop3}. \qed

Finally, we remark some holomorphy conditions of the system \eqref{eq:FS},\eqref{eq:FSH}.

{\bf Hamiltonian $H_{02}=\tilde{r}_{02}(H_{FS}-(p_1+p_2)), \ \tilde{r}_{02}:x=-((q_1-t)(p_1+p_2)-\alpha_2)(p_1+p_2), \ y=\frac{1}{p_1+p_2}, \ z=-((q_1-q_2)p_2-\alpha_0)p_2, \ w=\frac{1}{p_2}$}
\begin{align*}
&r_0^{02}:x_0=q_1, \ y_0=p_1, \ z_0=-(q_2 p_2-\alpha_0)p_2, \ w_0=\frac{1}{p_2}, \\
&r_1^{02}:x_1=\frac{1}{q_1}, \ y_1=-((p_1-p_2)q_1+\alpha_1)q_1, \ z_1=q_2+q_1, \ w_1=p_2, \\
&r_2^{02}:x_2=-(q_1p_1-\alpha_2)p_1, \ y_2=\frac{1}{p_1}, \ z_2=q_2 \ w_2=p_2, \\
&r_3^{02}:x_3=q_1+\frac{q_2 p_2+\alpha_3-(\alpha_0+\alpha_2+\eta)}{p_1}-\frac{t}{p_1^2}, \ y_3=p_1, \ z_3=q_2 p_1, \ w_3=\frac{p_2}{p_1}, \\
&r_4^{02}:x_4=q_1-q_2+\frac{2q_2p_2+\alpha_4-(\alpha_0+\alpha_2)}{p_1}+\frac{1-t}{p_1^2}, \ y_4=p_1, \ z_4=q_2 p_1^2, \ w_4=\frac{p_2-p_1}{p_1^2}, \\
&r_5^{02}:x_5=q_1, \ y_5=p_1, \ z_5=-(q_2 p_2-\alpha_5-\alpha_0)p_2, \ w_5=\frac{1}{p_2},\\
&r_6^{02}:x_6=-(q_1p_1+q_2 p_2-(\eta+\alpha_0+\alpha_2))p_1, \ y_6=\frac{1}{p_1}, \ z_6=q_2 p_1, \ w_6=\frac{p_2}{p_1},
\end{align*}
where $r_3^{02} \left(H_{02}+\frac{1}{p_1} \right), \ r_4^{02} \left(H_{02}+\frac{1}{p_1} \right)$ (cf. \eqref{reducedholo}).

\vspace{0.5cm}

{\bf Appendix A :Reformulation of Fuji-Suzuki coupled Painlev\'e VI system}

In this Appendix A, we will reformulate the  Hamiltonian system  \eqref{eq:FS},\eqref{eq:FSH} by replacing its constant complex parameters $\alpha_j \ (0,1,\ldots,5)$ and $\eta$ by $\beta_j \ (j=0,1,\ldots,6)$.

Define birational and symplectic transformations $r_i \ (i=0,1,\ldots,6)$ as follows:
\begin{align}
\begin{split}
r_0:(x_0,y_0,z_0,w_0)=&\left(-((q_1-q_2)p_1-\beta_0)p_1,\frac{1}{p_1},q_2,p_2+p_1 \right),\\
r_1:(x_1,y_1,z_1,w_1)=&\left(\frac{1}{q_1},-(q_1p_1+\beta_1)q_1,q_2,p_2 \right),\\
r_2:(x_2,y_2,z_2,w_2)=&\left(-((q_1-t)p_1-\beta_2)p_1,\frac{1}{p_1},q_2,p_2 \right),\\
r_3:(x_3,y_3,z_3,w_3)=&\left(-(q_1p_1+q_2p_2-\beta_3)p_1,\frac{1}{p_1},q_2p_1,\frac{p_2}{p_1} \right),\\
r_4:(x_4,y_4,z_4,w_4)=&\left(q_1,p_1,-((q_2-1)p_2-\beta_4)p_2,\frac{1}{p_2} \right),\\
r_5:(x_5,y_5,z_5,w_5)=&\left(q_1,p_1,\frac{1}{q_2},-(p_2q_2+\beta_5)q_2 \right),\\
r_6:(x_6,y_6,z_6,w_6)=&\left(-(XY+ZW-\beta_6)Y,\frac{1}{Y},ZY,\frac{W}{Y} \right),
\end{split}
\end{align}
where the coordinate system $(X,Y,Z,W)$ is given by
\begin{align*}
\begin{split}
&r_5 \circ r_1:(X,Y,Z,W):=\left(\frac{1}{q_1},-(q_1p_1+\alpha_1)q_1,\frac{1}{q_2},-(p_2q_2+\alpha_5)q_2 \right).
\end{split}
\end{align*}
There exist a polynomial $\tilde{H}_1$, such that the Hamiltonian system
\begin{equation}\label{eq:100}
   \frac{dq_1}{dt} =\frac{\partial \tilde{H}_1}{\partial p_1}, \ \frac{dp_1}{dt} =-\frac{\partial \tilde{H}_1}{\partial q_1},  \ \frac{dq_2}{dt} =\frac{\partial \tilde{H}_1}{\partial p_2}, \ \frac{dp_2}{dt} =-\frac{\partial \tilde{H}_1}{\partial q_2}
\end{equation}
is transformed into a polynomial Hamiltonian system under the action of each $r_i \ (i=0,1,\ldots,6)$, where a polynomial Hamiltonians $\tilde{H}_1$ is given by 
{\footnotesize
\begin{align}\label{eq:101}
\begin{split}
\tilde{H}_1=& \frac{H_{VI}(q_1,p_1;\beta_2,\beta_0+\beta_4,\beta_3+\beta_5,\beta_1(\beta_1+\beta_5+\beta_6))}{(\beta_0+2\beta_1+\beta_2+\beta_3+\beta_4+2\beta_5+\beta_6)t(t-1)}+\frac{H_{VI}(q_2,p_2;\beta_0+\beta_2,\beta_4,\beta_3+\beta_1,\beta_5(\beta_1+\beta_5+\beta_6))}{(\beta_0+2\beta_1+\beta_2+\beta_3+\beta_4+2\beta_5+\beta_6)t(t-1)}\\
&+\frac{(q_1-t)(q_2-1)\{(q_1p_1+\beta_1)p_2+p_1(q_2p_2+\beta_5)\}}{(\beta_0+2\beta_1+\beta_2+\beta_3+\beta_4+2\beta_5+\beta_6)t(t-1)},
\end{split}
\end{align}}
where $q_i,p_i \ (i=1,2)$ denote unknown complex variables, and $\beta_j \ (j=0,1,\ldots,6)$ are complex constant parameters satisfying the parameter's relation:
\begin{equation}\label{prarel11}
\beta_0+2\beta_1+\beta_2+\beta_3+\beta_4+2\beta_5+\beta_6=1.
\end{equation}

The relations between $\alpha_i \ (i=0,1,\ldots,5), \eta$ and $\beta_j \ (j=0,1,\ldots,6)$ are explicitly given as follows:
{\small
\begin{align}\label{prafss}
\begin{split}
&\alpha_0=\beta_0, \ \alpha_1=\beta_1, \ \alpha_2=\beta_2, \ \alpha_3=\beta_1+\beta_3+\beta_5+\beta_6, \ \alpha_4=\beta_4 , \alpha_5=\beta_5, \ \eta=\beta_1+\beta_5+\beta_6.
\end{split}
\end{align}}
Of course, $\alpha_i$ and $\beta_j$ satisfy the relation:
\begin{align*}
\begin{split}
&\alpha_0+\alpha_1+\alpha_2+\alpha_3+\alpha_4+\alpha_5=\beta_0+2\beta_1+\beta_2+\beta_3+\beta_4+2\beta_5+\beta_6=1.
\end{split}
\end{align*}

We remark that on new constant complex parameters $\beta_j \ (j=0,1,\ldots,6)$ the Hamiltonian system \eqref{eq:100},\eqref{eq:101} is invariant under these birational and symplectic transformations $s_0,s_1,\ldots,s_{10}$ (cf. Appendix B in  \cite{FS1}), whose generators are defined as follows$:$ with {\it the notation} $(*):=(q_1,p_1,q_2,p_2,t;\beta_0,\beta_1,\ldots,\beta_6)$;

{\footnotesize

\begin{align}\label{eq:103}
\begin{split}
s_0:(*) \rightarrow &\left(q_1,p_1-\frac{\beta_0}{q_1-q_2},q_2,p_2+\frac{\beta_0}{q_1-q_2},t;-\beta_0,\beta_1+\beta_0,\beta_2,\beta_3-\beta_0,\beta_4,\beta_5+\beta_0,\beta_6-\beta_0 \right),\\
s_1:(*) \rightarrow &\left(q_1+\frac{\beta_1}{p_1},p_1,q_2,p_2,t;\beta_0+\beta_1,-\beta_1,\beta_2+\beta_1,\beta_3+\beta_1,\beta_4,\beta_5,\beta_6+\beta_1 \right),\\
s_2:(*) \rightarrow &\left(q_1,p_1-\frac{\beta_2}{q_1-t},q_2,p_2,t;\beta_0,\beta_1+\beta_2,-\beta_2,\beta_3,\beta_4,\beta_5,\beta_6 \right),\\
s_3:(*) \rightarrow &(q_1+\frac{(\beta_1+\beta_3+\beta_5+\beta_6)q_1}{q_1p_1+q_2p_2-\beta_3},p_1-\frac{(\beta_1+\beta_3+\beta_5+\beta_6)p_1}{q_1p_1+q_2p_2+\beta_1+\beta_5+\beta_6},\\
&q_2+\frac{(\beta_1+\beta_3+\beta_5+\beta_6)q_2}{q_1p_1+q_2p_2-\beta_3},p_2-\frac{(\beta_1+\beta_3+\beta_5+\beta_6)p_2}{q_1p_1+q_2p_2+\beta_1+\beta_5+\beta_6},t;\beta_0,\beta_1,\\
&\beta_1+\beta_2+\beta_3+\beta_5+\beta_6,-\beta_1-\beta_5-\beta_6,\beta_1+\beta_3+\beta_4+\beta_5+\beta_6,\beta_5,-\beta_1-\beta_3-\beta_5),\\
s_4:(*) \rightarrow &\left(q_1,p_1,q_2,p_2-\frac{\beta_4}{q_2-1},t;\beta_0,\beta_1,\beta_2,\beta_3,-\beta_4,\beta_5+\beta_4,\beta_6 \right),\\
s_5:(*) \rightarrow &\left(q_1,p_1,q_2+\frac{\beta_5}{p_2},p_2,t;\beta_0+\beta_5,\beta_1,\beta_2,\beta_3+\beta_5,\beta_4+\beta_5,-\beta_5,\beta_6+\beta_5 \right),\\
s_6:(*) \rightarrow &\left(\frac{t}{q_2},-\frac{(p_2q_2+\beta_5)q_2}{t},\frac{t}{q_1},-\frac{(p_1q_1+\beta_1)q_1}{t},t;\beta_0,\beta_5,\beta_4,\beta_6,\beta_2,\beta_1,\beta_3 \right),\\
s_7:(*) \rightarrow &\left(\frac{1}{q_1},-(p_1q_1+\beta_1)q_1,\frac{1}{q_2},-(p_2q_2+\beta_5)q_2,\frac{1}{t};\beta_0,\beta_1,\beta_2,\beta_6,\beta_4,\beta_5,\beta_3 \right),\\
s_8:(*) \rightarrow &\left(\frac{q_2}{t},tp_2,\frac{q_1}{t},tp_1,\frac{1}{t};\beta_0,\beta_5,\beta_4,\beta_3,\beta_2,\beta_1,\beta_6 \right),\\
s_9:(*) \rightarrow &(\frac{1}{q_1},-(p_1q_1+p_2q_2+\beta_1+\beta_5+\beta_6)q_1,\frac{q_2}{q_1},p_2q_1,\frac{1}{t};\beta_4,\beta_1+\beta_3+\beta_5+\beta_6,\beta_2,\\
&-\beta_5-\beta_6,\beta_0,\beta_5,-\beta_3-\beta_5 ),\\
s_{10}:(*) \rightarrow &(\frac{t}{q_2},-\frac{(p_2q_2+p_1q_1+\beta_1+\beta_5+\beta_6)q_2}{t},\frac{q_1}{q_2},p_1q_2,t;\beta_2,\beta_1+\beta_3+\beta_5+\beta_6,\beta_4,\\
&-\beta_1-\beta_6,\beta_0,\beta_1,-\beta_1-\beta_3 ).
\end{split}
\end{align}}
We note that the subgroup $<s_0,s_1,\ldots,s_5>$ generated by $s_0,s_1,\ldots,s_5$ is isomorphic to the affine Weyl group of type $A_5^{(1)}$  (see Appendix B in  \cite{FS1}), and the transformation $s_6$ was found by Professor K. Fuji in Kobe university in August 2012.

We also remark that these transformations $s_i$ satisfies the following relations:
\begin{equation*}
s_{10}=s_8 \circ s_9, \quad s_6=s_8 \circ s_7, \quad s_k^2=1 \ (k=0,1,\ldots,9), \quad s_{10}^3=1.
\end{equation*}

Finally, let us define the following translation operators{\rm : \rm}
\begin{align*}
\begin{split}
&T_1:=(s_2s_{10}  s_{10}  s_1)^4, \quad T_2:=s_1T_1s_1,  \quad T_3:=s_5T_1s_5.
\end{split}
\end{align*}
These translation operators act on parameters $\beta_i$ as follows$:$
\begin{align*}
\begin{split}
T_1(\beta_0,\beta_1,\ldots,\beta_6)=&(\beta_0,\beta_1,\ldots,\beta_6)+(0,-1,1,0,-1,1,0),\\
T_2(\beta_0,\beta_1,\ldots,\beta_6)=&(\beta_0,\beta_1,\ldots,\beta_6)+(-1,1,0,-1,-1,1,-1),\\
T_3(\beta_0,\beta_1,\ldots,\beta_6)=&(\beta_0,\beta_1,\ldots,\beta_6)+(1,-1,1,1,0,-1,1).
\end{split}
\end{align*}

\vspace{0.5cm}

{\bf Appendix B :Searching for the B{\"a}cklund transformation $s_3$}

In this Appendix B, we will make Fuji-Suzuki's B{\"a}cklund transformation $s_3$ in \eqref{eq:103} by our method.

The key property is given as follows (see \cite{Sasa6,Sasa7}):
\begin{align}\label{moon2}
\begin{split}
&r:\left(-(XY+ZW-\beta),\frac{1}{Y},ZY,\frac{W}{Y} \right) \Longleftrightarrow s:\left(X,Y+\frac{ZW-\beta}{X},\frac{Z}{X},WX \right),\\
&r':\left(\frac{1}{X},-(XY+ZW+\beta)X,\frac{Z}{X},WX \right) \Longleftrightarrow s':\left(X+\frac{ZW+\beta}{Y},Y,ZY,\frac{W}{Y} \right).
\end{split}
\end{align}
These transformations $r,r',s$ and $s'$ are birational and symplectic, however, these are not auto-B{\"a}cklund transformations. These transformations can be considered as a relation between symmetry and holomorphy conditions appearing in the Garnier system (see \cite{Sasa6,Sasa7}).

\begin{tabular}{|c||c|} \hline
Equations & Relation between symmetry and holomorphy conditions  \\ \hline
Painlev\'e equations & \eqref{moon}  \\ \hline
Garnier systems  & \eqref{moon2}    \\ \hline
\end{tabular}

At first, for the system \eqref{eq:100},\eqref{eq:101}, we will make the above transformation.
\begin{proposition}
The birational and symplectic transformation
\begin{align}\label{eq:5-1}
\begin{split}
S:(q_1,p_1,q_2,p_2) \rightarrow \left(q_1,p_1+\frac{q_2 p_2-\beta_3}{q_1},\frac{q_2}{q_1},p_2 q_1 \right)
\end{split}
\end{align}
takes the system \eqref{eq:100},\eqref{eq:101} to a Hamiltonian system
\begin{equation}\label{eq:5-2}
   \frac{dq_1}{dt} =\frac{\partial H_2}{\partial p_1}, \ \frac{dp_1}{dt} =-\frac{\partial H_2}{\partial q_1},  \ \frac{dq_2}{dt} =\frac{\partial H_2}{\partial p_2}, \ \frac{dp_2}{dt} =-\frac{\partial H_2}{\partial q_2}
\end{equation}
with the polynomial Hamiltonian\rm{: \rm}
\begin{align}\label{eq:5-3}
\begin{split}
&(\beta_0+2\beta_1+\beta_2+\beta_3+\beta_4+2\beta_5+\beta_6)t(t-1) H_2=\\
&t p_2 p_1 - q_1 p_1 - p_2 q_1 p_1 + q_1^2 p_1 + t q_1 p_1^2 - q_1^2 p_1^2 - t q_1^2 p_1^2 + 
 q_1^3 p_1^2 + p_2^2 q_2 - p_2 q_1 q_2 - t p_2 p_1 q_2\\
&+ p_2 q_1 p_1 q_2 + t p_2 q_1 p_1 q_2 - p_2 q_1^2 p_1 q_2 -
  p_2^2 q_2^2 - t p_2^2 q_2^2 + p_2 q_1 q_2^2 - t p_2 q_1 p_1 q_2^2 + p_2 q_1^2 p_1 q_2^2 + 
 t p_2^2 q_2^3\\
& + t q_1 p_1 \beta_0 - q_1^2 p_1 \beta_0 + p_2 q_2 \beta_0 - 
 t p_2 q_2 \beta_0 + p_2 q_1 q_2 \beta_0 - p_2 q_1 q_2^2 \beta_0 - p_2 \beta_1 +
  p_2 q_2 \beta_1 + p_2 q_1 q_2 \beta_1\\
& - p_2 q_1 q_2^2 \beta_1 + 
 q_1 \beta_1^2 + q_1 p_1 \beta_2 - q_1^2 p_1 \beta_2 + p_2 q_1 q_2 \beta_2 - 
 p_2 q_1 q_2^2 \beta_2 - \beta_3 - p_2 \beta_3 + q_1 \beta_3 + 
 t p_1 \beta_3\\
& - q_1 p_1 \beta_3 - t q_1 p_1 \beta_3 + q_1^2 p_1 \beta_3 + 
 p_2 q_2 \beta_3 + t p_2 q_2 \beta_3 - t p_2 q_2^2 \beta_3 + 
 t \beta_0 \beta_3 - q_1 \beta_0 \beta_3 + \beta_2 \beta_3\\
& -
  q_1 \beta_2 \beta_3 + t q_1 p_1 \beta_4 - q_1^2 p_1 \beta_4 + 
 p_2 q_1 q_2 \beta_4 - p_2 q_1 q_2^2 \beta_4 + t \beta_3 \beta_4 - 
 q_1 \beta_3 \beta_4 + t q_1 p_1 \beta_5 - q_1^2 p_1 \beta_5\\
& - 
 t p_2 q_2 \beta_5 + p_2 q_1 q_2 \beta_5 - t q_1 p_1 q_2 \beta_5 + 
 q_1^2 p_1 q_2 \beta_5 + t p_2 q_2^2 \beta_5 - p_2 q_1 q_2^2 \beta_5 + 
 q_1 \beta_1 \beta_5 + q_1 q_2 \beta_1 \beta_5\\
& + 
 t \beta_3 \beta_5 - q_1 \beta_3 \beta_5 - 
 t q_2 \beta_3 \beta_5 + q_1 q_2 \beta_3 \beta_5 + q_1 q_2 \beta_5^2 +
  q_1 \beta_1 \beta_6 + q_1 q_2 \beta_5 \beta_6.
\end{split}
\end{align}
\end{proposition}
For this system \eqref{eq:5-2},\eqref{eq:5-3}, we can find the holomorphy condition:
\begin{align*}
R_3:(X_1,Y_1,Z_1,W_1)=&\left(\frac{1}{q_1},-(q_1p_1+\beta_1+\beta_3+\beta_5+\beta_6)q_1,q_2,p_2 \right).
\end{align*}
Next, we will explain how to find the B{\"a}cklund transformation $s_3$. Here, let us review the relation between symmetry and holomorphy (see \cite{YS}):
\begin{equation}\label{moon}
r:\left(\frac{1}{X},-(YX+\beta),Z,W \right) \Longleftrightarrow s:\left(X+\frac{\beta}{Y},Y,Z,W \right).
\end{equation}
By using this method, we can obtain the following B{\"a}cklund transformation for this system \eqref{eq:5-2},\eqref{eq:5-3}.

\begin{proposition}
The system \eqref{eq:5-2},\eqref{eq:5-3} is invariant under the following birational and symplectic transformation\rm{ : \rm}
\begin{align}
\begin{split}
S_3:(q_1,p_1,q_2,p_2,t;\beta_0,\beta_1,\ldots,\beta_6) \rightarrow &(q_1+\frac{\beta_1+\beta_3+\beta_5+\beta_6}{p_1},p_1,q_2,p_2,t;\beta_0,\beta_1,\\
&\beta_1+\beta_2+\beta_3+\beta_5+\beta_6,-\beta_1-\beta_5-\beta_6,\\
&\beta_1+\beta_3+\beta_4+\beta_5+\beta_6,\beta_5,-\beta_1-\beta_3-\beta_5).
\end{split}
\end{align}
\end{proposition}
Pulling back the transformation $S_3$ by the birational transformation \eqref{eq:5-1}, we can obtain Fuji-Suzuki's B{\"a}cklund transformation $s_3$ in \eqref{eq:103}.

\vspace{0.5cm}

{\bf Appendix C :Holomorphy conditions of type III}

\begin{figure}
\unitlength 0.1in
\begin{picture}( 56.6800, 43.8200)(  4.5000,-44.7200)
%
\special{pn 8}%
\special{pa 3136 778}%
\special{pa 1876 1226}%
\special{dt 0.045}%
\special{pa 1888 1220}%
\special{pa 4174 1220}%
\special{dt 0.045}%
\special{pa 3142 778}%
\special{pa 4162 1220}%
\special{dt 0.045}%
%
\special{pn 8}%
\special{pa 1902 1220}%
\special{pa 1494 2010}%
\special{dt 0.045}%
\special{pa 1494 2010}%
\special{pa 4644 2010}%
\special{dt 0.045}%
\special{pa 4162 1222}%
\special{pa 4644 2010}%
\special{dt 0.045}%
%
\special{pn 8}%
\special{pa 3142 786}%
\special{pa 3136 1494}%
\special{dt 0.045}%
\special{pa 1902 1216}%
\special{pa 2376 1686}%
\special{dt 0.045}%
\special{pa 2376 1686}%
\special{pa 3762 1686}%
\special{dt 0.045}%
\special{pa 4162 1222}%
\special{pa 3776 1688}%
\special{dt 0.045}%
\special{pa 2376 1688}%
\special{pa 3136 1494}%
\special{dt 0.045}%
\special{pa 3136 1494}%
\special{pa 3768 1682}%
\special{dt 0.045}%
\special{pa 1502 2010}%
\special{pa 2382 1688}%
\special{dt 0.045}%
\special{pa 3780 1686}%
\special{pa 4644 2010}%
\special{dt 0.045}%
%
\special{pn 8}%
\special{pa 1508 2010}%
\special{pa 4630 2010}%
\special{dt 0.045}%
%
\special{pn 20}%
\special{pa 3148 778}%
\special{pa 3312 738}%
\special{fp}%
\special{sh 1}%
\special{pa 3312 738}%
\special{pa 3242 734}%
\special{pa 3260 750}%
\special{pa 3252 772}%
\special{pa 3312 738}%
\special{fp}%
%
\special{pn 20}%
\special{pa 1502 2010}%
\special{pa 1444 2114}%
\special{fp}%
\special{sh 1}%
\special{pa 1444 2114}%
\special{pa 1494 2064}%
\special{pa 1470 2066}%
\special{pa 1460 2046}%
\special{pa 1444 2114}%
\special{fp}%
%
\special{pn 20}%
\special{pa 4644 2010}%
\special{pa 4688 2114}%
\special{fp}%
\special{sh 1}%
\special{pa 4688 2114}%
\special{pa 4680 2044}%
\special{pa 4666 2064}%
\special{pa 4642 2060}%
\special{pa 4688 2114}%
\special{fp}%
%
\special{pn 20}%
\special{pa 3142 1490}%
\special{pa 3326 1444}%
\special{fp}%
\special{sh 1}%
\special{pa 3326 1444}%
\special{pa 3256 1440}%
\special{pa 3274 1458}%
\special{pa 3266 1480}%
\special{pa 3326 1444}%
\special{fp}%
\put(33.5700,-7.6400){\makebox(0,0)[lb]{$\frac{1}{q_1}$}}%
\put(32.1200,-14.4000){\makebox(0,0)[lb]{$\frac{1}{q_2}$}}%
\put(46.9900,-21.4200){\makebox(0,0)[lb]{$\frac{1}{p_2}$}}%
\put(13.3600,-21.6600){\makebox(0,0)[lb]{$\frac{1}{p_1}$}}%
\put(29.8300,-11.0100){\makebox(0,0)[lb]{$r_6$}}%
%
\special{pn 8}%
\special{pa 3152 2796}%
\special{pa 1892 3244}%
\special{dt 0.045}%
\special{pa 1904 3238}%
\special{pa 4192 3238}%
\special{dt 0.045}%
\special{pa 3158 2796}%
\special{pa 4178 3238}%
\special{dt 0.045}%
%
\special{pn 8}%
\special{pa 1918 3238}%
\special{pa 1512 4030}%
\special{dt 0.045}%
\special{pa 1512 4030}%
\special{pa 4662 4030}%
\special{dt 0.045}%
\special{pa 4178 3242}%
\special{pa 4662 4030}%
\special{dt 0.045}%
%
\special{pn 8}%
\special{pa 3158 2806}%
\special{pa 3152 3512}%
\special{dt 0.045}%
\special{pa 1918 3236}%
\special{pa 2392 3704}%
\special{dt 0.045}%
\special{pa 2392 3704}%
\special{pa 3780 3704}%
\special{dt 0.045}%
\special{pa 4178 3242}%
\special{pa 3792 3708}%
\special{dt 0.045}%
\special{pa 2392 3708}%
\special{pa 3152 3512}%
\special{dt 0.045}%
\special{pa 3152 3512}%
\special{pa 3786 3700}%
\special{dt 0.045}%
\special{pa 1518 4030}%
\special{pa 2400 3708}%
\special{dt 0.045}%
\special{pa 3798 3704}%
\special{pa 4662 4030}%
\special{dt 0.045}%
%
\special{pn 8}%
\special{pa 1526 4030}%
\special{pa 4648 4030}%
\special{dt 0.045}%
%
\special{pn 20}%
\special{pa 3164 2796}%
\special{pa 3330 2758}%
\special{fp}%
\special{sh 1}%
\special{pa 3330 2758}%
\special{pa 3260 2754}%
\special{pa 3278 2770}%
\special{pa 3270 2792}%
\special{pa 3330 2758}%
\special{fp}%
%
\special{pn 20}%
\special{pa 1518 4030}%
\special{pa 1462 4134}%
\special{fp}%
\special{sh 1}%
\special{pa 1462 4134}%
\special{pa 1512 4084}%
\special{pa 1488 4086}%
\special{pa 1476 4066}%
\special{pa 1462 4134}%
\special{fp}%
%
\special{pn 20}%
\special{pa 4662 4030}%
\special{pa 4704 4134}%
\special{fp}%
\special{sh 1}%
\special{pa 4704 4134}%
\special{pa 4698 4064}%
\special{pa 4684 4084}%
\special{pa 4660 4080}%
\special{pa 4704 4134}%
\special{fp}%
%
\special{pn 20}%
\special{pa 3158 3508}%
\special{pa 3342 3464}%
\special{fp}%
\special{sh 1}%
\special{pa 3342 3464}%
\special{pa 3272 3460}%
\special{pa 3290 3476}%
\special{pa 3282 3498}%
\special{pa 3342 3464}%
\special{fp}%
\put(33.7400,-27.8300){\makebox(0,0)[lb]{$\frac{1}{q_1}$}}%
\put(32.2900,-34.5900){\makebox(0,0)[lb]{$\frac{1}{q_2}$}}%
\put(47.1600,-41.6200){\makebox(0,0)[lb]{$\frac{1}{p_2}$}}%
\put(13.5300,-41.8400){\makebox(0,0)[lb]{$\frac{1}{p_1}$}}%
\put(30.0000,-31.2100){\makebox(0,0)[lb]{$r_6$}}%
%
\special{pn 20}%
\special{pa 1508 1984}%
\special{pa 2360 1694}%
\special{fp}%
%
\special{pn 20}%
\special{pa 4174 1218}%
\special{pa 4640 2000}%
\special{fp}%
%
\special{pn 20}%
\special{pa 1792 1916}%
\special{pa 4390 1916}%
\special{fp}%
%
\special{pn 20}%
\special{pa 1604 1808}%
\special{pa 4502 1808}%
\special{fp}%
%
\special{pn 20}%
\special{pa 1914 1224}%
\special{pa 2350 1678}%
\special{fp}%
%
\special{pn 20}%
\special{pa 2094 1378}%
\special{pa 4020 1378}%
\special{fp}%
\put(42.4300,-16.9300){\makebox(0,0)[lb]{$r_5$}}%
\put(23.6800,-13.7000){\makebox(0,0)[lb]{$r_4$}}%
\put(19.7300,-15.7700){\makebox(0,0)[lb]{$r_1$}}%
\put(27.2000,-18.0000){\makebox(0,0)[lb]{$r_2$}}%
\put(31.9400,-20.3800){\makebox(0,0)[lb]{$r_0$}}%
\put(13.1900,-19.5400){\makebox(0,0)[lb]{$r_3$}}%
%
\special{pn 20}%
\special{pa 1914 3244}%
\special{pa 4174 3244}%
\special{fp}%
%
\special{pn 20}%
\special{pa 1526 4028}%
\special{pa 4648 4028}%
\special{fp}%
%
\special{pn 20}%
\special{pa 1612 3842}%
\special{pa 2256 3582}%
\special{fp}%
%
\special{pn 20}%
\special{pa 2402 3704}%
\special{pa 3788 3704}%
\special{fp}%
%
\special{pn 20}%
\special{pa 2730 3482}%
\special{pa 3168 4012}%
\special{fp}%
%
\special{pn 20}%
\special{pa 4012 3498}%
\special{pa 4296 3882}%
\special{fp}%
\put(23.3300,-32.3600){\makebox(0,0)[lb]{$r_1$}}%
\put(19.2100,-40.1200){\makebox(0,0)[lb]{$r_5$}}%
\put(17.6600,-37.2800){\makebox(0,0)[lb]{$r_0$}}%
\put(29.9600,-36.8100){\makebox(0,0)[lb]{$r_3$}}%
\put(39.7000,-37.8200){\makebox(0,0)[lb]{$r_4$}}%
\put(29.1000,-39.2700){\makebox(0,0)[lb]{$r_2$}}%
%
\special{pn 8}%
\special{pa 460 450}%
\special{pa 6118 450}%
\special{fp}%
\special{pa 890 464}%
\special{pa 890 464}%
\special{fp}%
%
\special{pn 8}%
\special{pa 450 2484}%
\special{pa 6110 2484}%
\special{fp}%
\special{pa 880 2500}%
\special{pa 880 2500}%
\special{fp}%
%
\special{pn 8}%
\special{pa 1216 158}%
\special{pa 1216 4472}%
\special{fp}%
%
\special{pn 8}%
\special{pa 5026 150}%
\special{pa 5026 4466}%
\special{fp}%
\put(20.6700,-4.1100){\makebox(0,0)[lb]{Accessible singular loci}}%
\put(6.7400,-4.1100){\makebox(0,0)[lb]{System}}%
\put(50.8000,-2.6000){\makebox(0,0)[lb]{Partition of}}%
\put(50.6100,-4.1800){\makebox(0,0)[lb]{Accessible singular loci}}%
\put(51.7200,-15.0800){\makebox(0,0)[lb]{$(a,b)=(3,3)$}}%
\put(51.8100,-34.6600){\makebox(0,0)[lb]{$(a,b)=(3,3)$}}%
\put(5.7900,-33.2000){\makebox(0,0)[lb]{$\eqref{th:1.4.2}$}}%
\put(5.8800,-15.0100){\makebox(0,0)[lb]{$\eqref{th:1.4.1}$}}%
\end{picture}%
\label{'SŽŸŒ³‹óŠÔFS}
\caption{\small This figure is the Hirzebruch manifold defined by H. Kimura (see \cite{K}). The bold lines denote some accessible singular loci for each system. Both systems are transformed by the birational transformation \eqref{Tr}. We also note that both types a,b are exchanged by the birational transformation \eqref{Tr}.}
\end{figure}

\begin{tabular}{|c||c|c|} \hline 
 & Holomorphy cond. of type a & Holomorphy cond. of type b \\ \hline
\eqref{th:1.4.1} & $r_1,r_3,r_5$ &  $r_0,r_2,r_4$ \\ \hline
\eqref{th:1.4.2} & $r_0,r_2,r_4$ &  $r_1,r_3,r_5$ \\ \hline
Type of Accessible sing. & ${\mathbb P}^1$ &  ${\mathbb P}^1$ \\ \hline
Type of Local index & $(2,1,0,1)$ &  $(2,1,1,0)$ \\ \hline
\end{tabular}

In this appendix, at first we will give some holomorphy conditions for the Hamiltonian system transformed the system \eqref{eq:FS},\eqref{eq:FSH} by the birational transformation;
\begin{equation}
sr_6:(Q_1,P_1,Q_2,P_2)=\left(q_1+\frac{q_2p_2+\eta}{p_1},p_1,\frac{p_2}{p_1},-q_2p_1 \right).
\end{equation}

Define birational and symplectic transformations $r_i \ (i=0,1,\ldots,6)$ as follows:
\begin{align}\label{th:1.4.1}
\begin{split}
r_0:(x_0,y_0,z_0,w_0)=&\left(-(Q_1P_1+(Q_2+1)P_2-\gamma_0)P_1,\frac{1}{P_1},(Q_2+1)P_1,\frac{P_2}{P_1} \right),\\
r_2:(x_2,y_2,z_2,w_2)=&\left(-((Q_1-t)P_1+Q_2P_2-\gamma_2)P_1,\frac{1}{P_1},Q_2P_1,\frac{P_2}{P_1} \right),\\
r_3:(x_3,y_3,z_3,w_3)=&\left(-(Q_1P_1-\gamma_3)P_1,\frac{1}{P_1},Q_2,P_2 \right),\\
r_5:(x_5,y_5,z_5,w_5)=&\left(Q_1,P_1,-(Q_2P_2-\gamma_5)P_2,\frac{1}{P_2} \right),\\
r_6:(x_6,y_6,z_6,w_6)=&\left(\frac{1}{Q_1},-(Q_1P_1+Q_2P_2+\gamma_6)Q_1,\frac{Q_2}{Q_1},P_2Q_1 \right),\\
r_1:(x_1,y_1,z_1,w_1)=&\left(-(x_6y_6-\gamma_1)y_6,\frac{1}{y_6},z_6,w_6 \right),\\
r_4:(x_4,y_4,z_4,w_4)=&\left(-(x_6y_6+(z_6-1)w_6-\gamma_4)y_6,\frac{1}{y_6},(z_6-1)y_6,\frac{w_6}{y_6} \right).
\end{split}
\end{align}
There exist a polynomial $\tilde{H}_1$, such that the Hamiltonian system
\begin{equation}\label{sys1}
   \frac{dQ_1}{dt} =\frac{\partial \tilde{H}_1}{\partial P_1}, \ \frac{dP_1}{dt} =-\frac{\partial \tilde{H}_1}{\partial Q_1},  \ \frac{dQ_2}{dt} =\frac{\partial \tilde{H}_1}{\partial P_2}, \ \frac{dP_2}{dt} =-\frac{\partial \tilde{H}_1}{\partial Q_2}
\end{equation}
is transformed into the polynomial Hamiltonian $\tilde{H}_1(Q_1,P_1,Q_2,P_2)=sr_6(H_{FS}(q_1,p_1,q_2,p_2))$.

The relations between $\alpha_i \ (i=0,1,\ldots,5), \eta$ and $\gamma_j \ (j=0,1,\ldots,6)$ are explicitly given as follows:
\begin{align}\label{pra3}
\begin{split}
&\gamma_0=\alpha_0+\eta, \quad \gamma_1=\alpha_1, \quad \gamma_2=\alpha_2+\eta, \quad \gamma_3=\alpha_3\\
&\gamma_4=\alpha_4+\eta, \quad \gamma_5=\alpha_5, \quad \gamma_6=-\eta.
\end{split}
\end{align}

Next, we will give some holomorphy conditions for the Hamiltonian system transformed the system \eqref{eq:FS},\eqref{eq:FSH} by the birational transformation;
\begin{equation}
rr_6:(\tilde{Q}_1,\tilde{P}_1,\tilde{Q}_2,\tilde{P}_2)=\left( \frac{q_1}{q_2},p_1q_2,\frac{1}{q_2},-(q_2p_2+q_1p_1+\eta)q_2 \right).
\end{equation}

Define birational and symplectic transformations $r_i \ (i=0,1,\ldots,6)$ as follows:
\begin{align}\label{th:1.4.2}
\begin{split}
r_0:(x_0,y_0,z_0,w_0)=&\left(-((\tilde{Q}_1-1)	\tilde{P}_1-(\gamma_0+\gamma_6))	\tilde{P}_1,\frac{1}{	\tilde{P}_1},\tilde{Q}_2,\tilde{P}_2 \right),\\
r_2:(x_2,y_2,z_2,w_2)=&\left(-((\tilde{Q}_1-t \tilde{Q}_2)	\tilde{P}_1-(\gamma_2+\gamma_6))	\tilde{P}_1,\frac{1}{	\tilde{P}_1},\tilde{Q}_2,\tilde{P}_2+t	\tilde{P}_1 \right),\\
r_4:(x_4,y_4,z_4,w_4)=&\left(\tilde{Q}_1,	\tilde{P}_1,-((\tilde{Q}_2-1)\tilde{P}_2-(\gamma_4+\gamma_6))\tilde{P}_2,\frac{1}{\tilde{P}_2} \right),\\
r_5:(x_5,y_5,z_5,w_5)=&\left(-(\tilde{Q}_1	\tilde{P}_1+\tilde{Q}_2\tilde{P}_2-(\gamma_5+\gamma_6))	\tilde{P}_1,\frac{1}{	\tilde{P}_1},\tilde{Q}_2\tilde{P}_1,\frac{\tilde{P}_2}{\tilde{P}_1} \right),\\
r_6:(x_6,y_6,z_6,w_6)=&\left(\frac{\tilde{Q}_1}{\tilde{Q}_2},\tilde{P}_1\tilde{Q}_2,\frac{1}{\tilde{Q}_2},-(\tilde{Q}_2\tilde{P}_2+\tilde{Q}_1\tilde{P}_1-\gamma_6)\tilde{Q}_2 \right),\\
r_1:(x_1,y_1,z_1,w_1)=&\left(-(XY+ZW-(\gamma_1+\gamma_6))Y,\frac{1}{Y},ZY,\frac{W}{Y} \right),\\
r_3:(x_3,y_3,z_3,w_3)=&\left(-(x_6y_6+z_6w_6-(\gamma_3+\gamma_6))y_6,\frac{1}{y_6},z_6y_6,\frac{w_6}{y_6} \right),
\end{split}
\end{align}
where the coordinate system $(X,Y,Z,W)$ is given by
\begin{align*}
\begin{split}
&R_6:(X,Y,Z,W)=\left(\frac{1}{\tilde{Q}_1},-(\tilde{Q}_1\tilde{P}_1+\tilde{Q}_2\tilde{P}_2-\gamma_6)\tilde{Q}_1,\frac{\tilde{Q}_2}{\tilde{Q}_1},\tilde{P}_2\tilde{Q}_1 \right).
\end{split}
\end{align*}
There exist a polynomial $\tilde{H}_6$, such that the Hamiltonian system
\begin{equation}\label{sys11}
   \frac{d \tilde{Q}_1}{dt} =\frac{\partial \tilde{H}_6}{\partial \tilde{P}_1}, \ \frac{d\tilde{P}_1}{dt} =-\frac{\partial \tilde{H}_6}{\partial \tilde{Q}_1},  \ \frac{d\tilde{Q}_2}{dt} =\frac{\partial \tilde{H}_6}{\partial \tilde{P}_2}, \ \frac{d\tilde{P}_2}{dt} =-\frac{\partial \tilde{H}_6}{\partial \tilde{Q}_2}
\end{equation}
is transformed into the polynomial Hamiltonian $\tilde{H}_6(\tilde{Q}_1,\tilde{P}_1,\tilde{Q}_2,\tilde{P}_2)=rr_6(H_{FS}(q_1,p_1,q_2,p_2))$ with parameter relations \eqref{pra3}.

We note that the condition $r_2$ should be read that $r_2(K-\tilde{P}_1\tilde{Q}_2)$ is a polynomial with respect to $x_2,y_2,z_2,w_2$.

We show that both systems \eqref{sys1},\eqref{sys11} are transformed by the birational transformation;
\begin{equation}\label{Tr}
Tr:(\tilde{Q}_1,\tilde{P}_1,\tilde{Q}_2,\tilde{P}_2)=\left(-\left(Q_2+\frac{Q_1P_1+\gamma_6}{P_2} \right),-P_2,-\frac{P_1}{P_2},Q_1P_2 \right).
\end{equation}

\begin{tabular}{|c||c|c|} \hline 
Holomorphy conditions & Parameter of \eqref{th:1.4.1} & Parameter of \eqref{th:1.4.2} \\ \hline
$r_0$ & $\gamma_0$ &  $\gamma_0+\gamma_6$ \\ \hline
$r_1$ & $\gamma_1$ &  $\gamma_1+\gamma_6$ \\ \hline
$r_2$ & $\gamma_2$ &  $\gamma_2+\gamma_6$ \\ \hline
$r_3$ & $\gamma_3$ &  $\gamma_3+\gamma_6$ \\ \hline
$r_4$ & $\gamma_4$ &  $\gamma_4+\gamma_6$ \\ \hline
$r_5$ & $\gamma_5$ &  $\gamma_5+\gamma_6$ \\ \hline
$r_6$ & $\gamma_6$ &  $-\gamma_6$ \\ \hline
\end{tabular}

\section{Appendix D}

{\bf Holomorphy conditions}

Define birational and symplectic transformations $r_i \ (i=0,1,\ldots,6)$ as follows:

\begin{align}\label{holo;FS2}
\begin{split}
r_0:(x_0,y_0,z_0,w_0)=&\left(-((q_1-q_2)p_1-\beta_0)p_1,\frac{1}{p_1},q_2,p_2+p_1 \right),\\
r_1:(x_1,y_1,z_1,w_1)=&\left(\frac{1}{q_1},-(q_1p_1+\beta_1)q_1,q_2,p_2 \right),\\
r_2:(x_2,y_2,z_2,w_2)=&\left(-((q_1-t)p_1-\beta_2)p_1,\frac{1}{p_1},q_2,p_2 \right),\\
r_3:(x_3,y_3,z_3,w_3)=&\left(-(q_1p_1+q_2p_2-\beta_3)p_1,\frac{1}{p_1},q_2p_1,\frac{p_2}{p_1} \right),\\
r_4:(x_4,y_4,z_4,w_4)=&\left(q_1,p_1,-((q_2-s)p_2-\beta_4)p_2,\frac{1}{p_2} \right),\\
r_5:(x_5,y_5,z_5,w_5)=&\left(q_1,p_1,\frac{1}{q_2},-(p_2q_2+\beta_5)q_2 \right),\\
r_6:(x_6,y_6,z_6,w_6)=&\left(-((q_1-1)p_1+(q_2-1)p_2-\beta_6)p_1,\frac{1}{p_1},(q_2-1)p_1,\frac{p_2}{p_1} \right).
\end{split}
\end{align}
There exist two polynomials $H_1^{FS}$ and $H_2^{FS}$, such that the Hamiltonian system
\begin{equation}\label{eq:FS2}
  \left\{
  \begin{aligned}
   dq_1 =&\frac{\partial H_1^{FS}}{\partial p_1}dt+\frac{\partial H_2^{FS}}{\partial p_1}ds, \quad dp_1 =-\frac{\partial H_1^{FS}}{\partial q_1}dt-\frac{\partial H_2^{FS}}{\partial q_1}ds,\\
   dq_2 =&\frac{\partial H_1^{FS}}{\partial p_2}dt+\frac{\partial H_2^{FS}}{\partial p_2}ds, \quad dp_2 =-\frac{\partial H_1^{FS}}{\partial q_2}dt-\frac{\partial H_2^{FS}}{\partial q_2}ds
   \end{aligned}
  \right.
\end{equation}
is transformed into a polynomial Hamiltonian system under the action of each $r_i \ (i=0,1,\ldots,6)$, where two polynomial Hamiltonians $H_1^{FS},H_2^{FS}$ are given by {\rm (cf. \cite{FS1,KNS1,KNS2}) \rm}
\begin{align}
\begin{split}
H_1^{FS} &=H_{VI}(q_1,p_1,t,s;\beta_2,\beta_0+\beta_4,\beta_1,\beta_5+\beta_6,\beta_3+\beta_5)\\
&+H_{VI}(q_2,p_2,t,s;\beta_0+\beta_2,\beta_4,\beta_5,\beta_1+\beta_6,\beta_1+\beta_3)\\
&+\frac{(q_1-t)(q_2-s) \{2(q_1 p_1+\beta_1)(q_2p_2+\beta_5)-(q_1 p_1+\beta_1)p_2-(q_2p_2+\beta_5)p_1 \}}{(\beta_0+2\beta_1+\beta_2+\beta_3+\beta_4+2\beta_5+\beta_6)t(t-1)(t-s)}\\
&-\frac{ 2\beta_1 \beta_5 s }{(\beta_0+2\beta_1+\beta_2+\beta_3+\beta_4+2\beta_5+\beta_6)(t-1)(t-s)},\\
H_2^{FS} &=\pi(H_1^{FS}), \quad \pi=\{q_1 \leftrightarrow q_2, \ p_1 \leftrightarrow p_2, \  t \leftrightarrow s, \ \beta_1 \leftrightarrow \beta_5, \ \beta_2 \leftrightarrow \beta_4 \}.
\end{split}
\end{align}
The symbol $H_{VI}(q,p,t,\eta;\gamma_0,\gamma_1,\gamma_2,\gamma_3,\gamma_4)$ denotes (see \cite{Sasa5})
\begin{align}
\begin{split}
&t(t-1)(t-\eta) H_{VI}(q,p,t,\eta;\gamma_0,\gamma_1,\gamma_2,\gamma_3,\gamma_4)\\
&=q(q-1)(q-\eta)(q-t)p^2+\{\gamma_1(t-\eta)q(q-1)+2\gamma_2 q(q-1)(q-\eta)\\
&+\gamma_3 (t-1)q(q-\eta)+\gamma_4 t(q-1)(q-\eta) \}p\\
&+\gamma_2 \{(\gamma_1+\gamma_2)(t-\eta)+\gamma_2(q-1)+\gamma_3(t-1)+t \gamma_4 \}q, \quad (\gamma_0+\gamma_1+2\gamma_2+\gamma_3+\gamma_4=1).
\end{split}
\end{align}

\vspace{0.5cm}

We note that the holomorphy conditions should be read that in the Hamiltonian $H_1^{FS}$
\begin{align*}
\begin{split}
&r_2(H_1^{FS} - p_1)
\end{split}
\end{align*}
are polynomials with respect to $x_2,y_2,z_2,w_2$, and in the Hamiltonian $H_2^{FS}$
\begin{align*}
\begin{split}
&r_4(H_2^{FS} -p_2)
\end{split}
\end{align*}
are polynomials with respect to $x_4,y_4,z_4,w_4$.

We see that the birational and symplectic transformation $\varphi$\rm{:\rm}
\begin{equation}
  \left\{
  \begin{aligned}
   Q_1=&\frac{1-q_1}{q_1}, \quad P_1=-(p_1 q_1+\beta_1)q_1, \quad Q_2=\frac{1-q_2}{q_2}, \quad P_2=-(p_2 q_2+\beta_5)q_2,\\
   T=&\frac{1-t}{t}, \quad S=\frac{1-s}{s}
   \end{aligned}
  \right. 
\end{equation}
takes the system \eqref{eq:FS2} into Fuji-Suzuki system (see \cite{FS1}) when $S=1$.

We remark that the relations between $\alpha_i \ (i=0,1,\ldots,5), \eta$ (see \cite{FS1}) and $\beta_j \ (j=0,1,\ldots,6)$ are explicitly given as follows:
\begin{align}\label{FSprelation}
\begin{split}
&\alpha_3=\beta_1+\beta_3+\beta_5+\beta_6, \quad \eta=\beta_1+\beta_3+\beta_5,\\
&\alpha_0=\beta_0, \ \alpha_1=\beta_1, \ \alpha_2=\beta_2, \ \alpha_4=\beta_4 , \alpha_5=\beta_5.
\end{split}
\end{align}
Of course, $\alpha_i$ and $\beta_j$ satisfy the relation:
\begin{align}
\begin{split}
&\alpha_0+\alpha_1+\alpha_2+\alpha_3+\alpha_4+\alpha_5=\beta_0+2\beta_1+\beta_2+\beta_3+\beta_4+2\beta_5+\beta_6=1.
\end{split}
\end{align}

\vspace{0.5cm}

{\bf Completely integrable}

{\footnotesize
\begin{proposition}
Setting
\begin{equation}
K_1:=-H_1+\frac{\beta_3(\beta_1+\beta_5)( \rm{Log \rm}(s-t)-\rm{Log \rm}(s-1) )}{(\beta_0+2\beta_1+\beta_2+\beta_3+\beta_4+2\beta_5+\beta_6)(t-1)^2}, \quad K_2:=-H_2.
\end{equation}
Two Hamiltonians $K_1$ and $K_2$ satisfy
\begin{equation}
\{K_1,K_2\}+\left(\frac{\partial}{\partial s} \right)K_1-\left(\frac{\partial}{\partial t} \right)K_2=0,
\end{equation}
where  $\{,\}$ denotes the Poisson brackets:
\begin{equation}
\{L_1,L_2\}=\frac{\partial L_1}{\partial p_1}\frac{\partial L_2}{\partial q_1}-\frac{\partial L_1}{\partial q_1}\frac{\partial L_2}{\partial p_1}+\frac{\partial L_1}{\partial p_2}\frac{\partial L_2}{\partial q_2}-\frac{\partial L_1}{\partial q_2}\frac{\partial L_2}{\partial p_2}.
\end{equation}
\end{proposition}
}

We remark that on new constant complex parameters $\beta_j \ (j=0,1,\ldots,6)$ the Hamiltonian system \eqref{eq:FS2} is invariant under these birational and symplectic transformations $s_0,s_1,\ldots,s_9$ (cf. Appendix B in  \cite{FS1}), whose generators are defined as follows$:$ with {\it the notation} $(*):=(q_1,p_1,q_2,p_2,t,s;\beta_0,\beta_1,\ldots,\beta_6)$;

{\footnotesize

\begin{align}\label{symmetry:FS}
\begin{split}
&s_0:(*) \rightarrow \left(q_1,p_1-\frac{\beta_0}{q_1-q_2},q_2,p_2+\frac{\beta_0}{q_1-q_2},t,s;-\beta_0,\beta_1+\beta_0,\beta_2,\beta_3-\beta_0,\beta_4,\beta_5+\beta_0,\beta_6-\beta_0 \right),\\
&s_1:(*) \rightarrow \left(q_1+\frac{\beta_1}{p_1},p_1,q_2,p_2,t,s;\beta_0+\beta_1,-\beta_1,\beta_2+\beta_1,\beta_3+\beta_1,\beta_4,\beta_5,\beta_6+\beta_1 \right),\\
&s_2:(*) \rightarrow \left(q_1,p_1-\frac{\beta_2}{q_1-t},q_2,p_2,t,s;\beta_0,\beta_1+\beta_2,-\beta_2,\beta_3,\beta_4,\beta_5,\beta_6 \right),\\
&s_3:(*) \rightarrow (\frac{q_1 (-q_1 p_1 + q_1^2 p_1 - p_2 q_2 + p_2 q_2^2 - \beta_1 + 
    q_1 \beta_1 - \beta_5 + q_2 \beta_5 - \beta_6)}{g_1},\\
& -\frac{g_1 (q_1 p_1^2 -
        q_1^2 p_1^2 + p_2 p_1 q_2 - p_2 p_1 q_2^2 + \beta_1^2 - p_1 \beta_3 + 
       q_1 p_1 \beta_3 + \beta_1 \beta_3 + q_1 p_1 \beta_5 - 
       p_1 q_2 \beta_5 + \beta_1 \beta_5 + 
       q_1 p_1 \beta_6 + \beta_1 \beta_6)}{(-q_1 p_1 + q_1^2 p_1 - p_2 q_2 + 
      p_2 q_2^2 + q_1 \beta_1 + \beta_3 + q_2 \beta_5) (-q_1 p_1 + q_1^2 p_1 - 
      p_2 q_2 + p_2 q_2^2 - \beta_1 + q_1 \beta_1 - \beta_5 + 
      q_2 \beta_5 - \beta_6)},\\
& \frac{
 q_2 (-q_1 p_1 + q_1^2 p_1 - p_2 q_2 + p_2 q_2^2 - \beta_1 + q_1 \beta_1 - \beta_5 +
     q_2 \beta_5 - \beta_6)}{g_2},\\
& -\frac{g_2 (p_2 q_1 p_1 - 
       p_2 q_1^2 p_1 + p_2^2 q_2 - p_2^2 q_2^2 - p_2 q_1 \beta_1 + p_2 q_2 \beta_1 - 
       p_2 \beta_3 + 
       p_2 q_2 \beta_3 + \beta_1 \beta_5 + \beta_3 \beta_5 + \
\beta_5^2 + p_2 q_2 \beta_6 + \beta_5 \beta_6)}{(-q_1 p_1 + q_1^2 p_1 - 
      p_2 q_2 + p_2 q_2^2 + q_1 \beta_1 + \beta_3 + q_2 \beta_5) (-q_1 p_1 + 
      q_1^2 p_1 - p_2 q_2 + p_2 q_2^2 - \beta_1 + q_1 \beta_1 - \beta_5 + 
      q_2 \beta_5 - \beta_6)},\\
&t,s;\beta_0,\beta_1,\beta_1+\beta_2+\beta_3+\beta_5+\beta_6,-\beta_1-\beta_5-\beta_6,\beta_1+\beta_3+\beta_4+\beta_5+\beta_6,\beta_5,-\beta_1-\beta_3-\beta_5),\\
&s_4:(*) \rightarrow \left(q_1,p_1,q_2,p_2-\frac{\beta_4}{q_2-s},t,s;\beta_0,\beta_1,\beta_2,\beta_3,-\beta_4,\beta_5+\beta_4,\beta_6 \right),\\
&s_5:(*) \rightarrow \left(q_1,p_1,q_2+\frac{\beta_5}{p_2},p_2,t,s;\beta_0+\beta_5,\beta_1,\beta_2,\beta_3+\beta_5,\beta_4+\beta_5,-\beta_5,\beta_6+\beta_5 \right),\\
&s_6:(*) \rightarrow \left(1-q_1,-p_1,1-q_2,-p_2,1-t,1-s;\beta_0,\beta_1,\beta_2,\beta_6,\beta_4,\beta_5,\beta_3 \right),\\
&s_7:(*) \rightarrow \left(q_2,p_2,q_1,p_1,s,t;\beta_0,\beta_5,\beta_4,\beta_3,\beta_2,\beta_1,\beta_6 \right),\\
&s_8:(*) \rightarrow (1 - q_1, \frac{q_1 p_1 - q_1^2 p_1 + p_2 q_2 - p_2 q_2^2 - \beta_3 + q_1 \beta_3 + 
 q_1 \beta_5 - q_2 \beta_5 + 
 q_1 \beta_6}{(-1 + q_1) q_1}, \frac{(-1 + s) (-1 + q_1) q_2}{
s q_1 + q_2 - s q_2 - q_1 q_2}, \\
&-\frac{(s q_1 + q_2 - s q_2 - q_1 q_2) (-s p_2 q_1 - p_2 q_2 + s p_2 q_2 + 
    p_2 q_1 q_2 - \beta_5 + s \beta_5 + q_1 \beta_5)}{(-1 + s) s (-1 + 
    q_1) q_1},\\
&1-t,1-s;\beta_4,\beta_1+\beta_3+\beta_5+\beta_6,\beta_2,-\beta_5-\beta_6,\beta_0,\beta_5,-\beta_3-\beta_5 ),\\
&s_9:=s_7 \circ s_8 \quad ((s_{9})^6=1),
\end{split}
\end{align}
where 
\begin{align*}
\begin{split}
&g_1:=-q_1 p_1 + 
  q_1^2 p_1 - p_2 q_2 + p_2 q_2^2 + \beta_3 - q_1 \beta_3 - q_1 \beta_5 + 
  q_2 \beta_5 - 
  q_1 \beta_6,\\
&g_2:=-q_1 p_1 + q_1^2 p_1 - p_2 q_2 + p_2 q_2^2 + 
  q_1 \beta_1 - q_2 \beta_1 + \beta_3 - q_2 \beta_3 - 
  q_2 \beta_6.
\end{split}
\end{align*}
We note that the subgroup $<s_0,s_1,\ldots,s_5>$ generated by $s_0,s_1,\ldots,s_5$ is isomorphic to the affine Weyl group of type $A_5^{(1)}$  (see Appendix B in  \cite{FS1}), and the transformation $s_6$ was found by Professor K. Fuji in Kobe university in August 2012.

Finally, let us define the following translation operators{\rm : \rm}
\begin{align}
\begin{split}
&T_1:=(s_2s_9  s_9  s_1)^4, \quad T_2:=s_1T_1s_1,  \quad T_3:=s_5T_1s_5.
\end{split}
\end{align}
These translation operators act on parameters $\beta_i$ as follows$:$
\begin{align}
\begin{split}
T_1(\beta_0,\beta_1,\ldots,\beta_6)=&(\beta_0,\beta_1,\ldots,\beta_6)+(0,-1,1,0,-1,1,0),\\
T_2(\beta_0,\beta_1,\ldots,\beta_6)=&(\beta_0,\beta_1,\ldots,\beta_6)+(-1,1,0,-1,-1,1,-1),\\
T_3(\beta_0,\beta_1,\ldots,\beta_6)=&(\beta_0,\beta_1,\ldots,\beta_6)+(1,-1,1,1,0,-1,1).
\end{split}
\end{align}}

Finally, we remark some holomorphy conditions of the system \eqref{eq:FS2}.

{\bf Hamiltonians $H_{04}^{(1)}=r_{04}(H_1^{FS}), \ H_{04}^{(2)}=r_{04}(H_2^{FS}-p_1-p_2), \ r_{04}:x=-((q_1-q_2)p_1-\beta_0)p_1, \ y=\frac{1}{p_1}, \ z=-((q_2-s)(p_2+p_1)-\beta_4)(p_2+p_1), \ w=\frac{1}{p_2+p_1}$}
\begin{align*}
&r_0^{04}:x_0=-(q_1p_1-\beta_0)p_1, \ y_0=\frac{1}{p_1}, \ z_0=q_2, \ w_0=p_2, \\
&r_1^{04}:x_1=-(q_1p_1-\beta_1-\beta_0)p_1, \ y_1=\frac{1}{p_1}, \ z_1=q_2, \ w_1=p_2, \\
&r_2^{04}:x_2=q_1-q_2+\frac{2q_2p_2+\beta_2-(\beta_0+\beta_4)}{p_1}+\frac{t-s}{p_1^2}, \ y_2=p_1, \ z_2=q_2 p_1^2 \ w_2=\frac{p_2-p_1}{p_1^2}, \\
&r_3^{04}:x_3=q_1p_2, \ y_3=\frac{p_1}{p_2}, \ z_3=q_2+\frac{q_1p_1+\beta_3-(\beta_0+\beta_4)}{p_2}-\frac{s}{p_2^2}, \ w_3=p_2, \\
&r_4^{04}:x_4=q_1, \ y_4=p_1, \ z_4=-(q_2 p_2-\beta_4)p_2, \ w_4=\frac{1}{p_2}, \\
&r_5^{04}:x_5=\frac{1}{q_1}, \ y_5=-((p_1-p_2)q_1+\beta_5)q_1, \ z_5=q_2+q_1, \ w_5=p_2,\\
&r_6^{04}:x_6=q_1p_2, \ y_6=\frac{p_1}{p_2}, \ z_6=q_2+\frac{q_1p_1+\beta_6-(\beta_0+\beta_4)}{p_2}-\frac{s-1}{p_2^2}, \ w_6=p_2,
\end{align*}
where $r_2^{04} \left(H_{04}^{(1)}-\frac{1}{p_1} \right), \ r_2^{04} \left(H_{04}^{(2)}+\frac{1}{p_1} \right), \ r_3^{04} \left(H_{04}^{(2)}+\frac{1}{p_2} \right), \ r_6^{04} \left(H_{04}^{(2)}+\frac{1}{p_2} \right)$.

\hspace{0.5cm}

{\bf Hamiltonians $H_{045}^{(1)}=r_5^{04}(H_{04}^{(1)}), \ H_{045}^{(2)}=r_5^{04}(H_{04}^{(2)})$}

\begin{align*}
&r_0^{045}:x_0=\frac{1}{q_1}, \ y_0=-(p_1 q_1+\beta_0+\beta_5)q_1, \ z_0=q_2, \ w_0=p_2, \\
&r_1^{045}:x_1=\frac{1}{q_1}, \ y_1=-(p_1 q_1+\beta_1+\beta_0+\beta_5)q_1, \ z_1=q_2, \ w_1=p_2, \\
&r_2^{045}:x_2=x_5 w_5^2, \ y_2=\frac{y_5}{w_5^2}, \ z_2=z_5+\frac{2x_5 y_5+\beta_2-(\beta_0+\beta_4)}{w_5}+\frac{t-s}{w_5^2} \ w_2=w_5, \\
&r_3^{045}:x_3=\frac{q_1}{p_2}, \ y_3=p_1 p_2, \ z_3=q_2-\frac{q_1p_1-\beta_3+\beta_0+\beta_4+\beta_5}{p_2}-\frac{s}{p_2^2}, \ w_3=p_2, \\
&r_4^{045}:x_4=q_1, \ y_4=p_1-\frac{\beta_4 q_2}{q_1q_2-1}, \ z_4=q_2, \ w_4=p_2-\frac{\beta_4 q_1}{q_1q_2-1}, \\
&r_5^{045}:x_5=\frac{1}{q_1}, \ y_5=-(p_1 q_1+\beta_5)q_1, \ z_5=q_2, \ w_5=p_2,\\
&r_6^{045}:x_6=\frac{q_1}{p_2}, \ y_6=p_1 p_2, \ z_6=q_2-\frac{q_1p_1-\beta_6+\beta_0+\beta_4+\beta_5}{p_2}-\frac{s-1}{p_2^2}, \ w_6=p_2,
\end{align*}
where $r_2^{045} \left(r_5^{045}(H_{045}^{(1)})-\frac{1}{w_5} \right), \ r_2^{045} \left(r_5^{045}(H_{045}^{(2)})+\frac{1}{w_5} \right), \ r_3^{045} \left(H_{045}^{(2)}+\frac{1}{p_2} \right), \ r_6^{045} \left(H_{045}^{(2)}+\frac{1}{p_2} \right)$.

{\bf Hamiltonians $H_{15}^{(1)}=r_{01}\left(-\frac{1}{T^2} H_1^{FS} \right), \ H_{15}^{(2)}=r_{15}\left(-\frac{1}{S^2} H_2^{FS} \right), \ r_{15}:Q_1=\frac{1-q_1}{q_1}, \ P_1=-(p_1 q_1+\beta_1)q_1, \ Q_2=\frac{1-q_2}{q_2}, \ P_2=-(p_2 q_2+\beta_5)q_2, \ T=\frac{1-t}{t}, \ S=\frac{1-s}{s}$}
\begin{align*}
&r_0^{15}:x_0=-((q_1-q_2)p_1-\beta_0)p_1, \ y_0=\frac{1}{p_1}, \ z_0=q_2, \ w_0=p_2+p_1, \\
&r_1^{15}:x_1=\frac{1}{q_1}, \ y_1=-(q_1p_1+\beta_1)q_1, \ z_1=q_2, \ w_1=p_2, \\
&r_2^{15}:x_2=-((q_1-t)p_1-\beta_2)p_1, \ y_2=\frac{1}{p_1}, \ z_2=q_2, \ w_2=p_2, \\
&r_3^{15}:x_3=\frac{1}{q_1}, \ y_3=-(p_1 q_1+p_2 q_2+\beta_3+\beta_1+\beta_5)q_1, \ z_3=\frac{q_2}{q_1}, \ w_3=p_2 q_1, \\
&r_4^{15}:x_4=q_1, \ y_4=p_1, \ z_4=-((q_2-s)p_2-\beta_4)p_2, \ w_4=\frac{1}{p_2}, \\
&r_5^{15}:x_5=q_1, \ y_5=p_1, \ z_5=\frac{1}{q_2}, \ w_5=-(p_2q_2+\beta_5)q_2,\\
&r_6^{15}:x_6=-(q_1 p_1+q_2 p_2-\beta_6)p_1, \ y_6=\frac{1}{p_1}, \ z_6=q_2 p_1, \ w_6=\frac{p_2}{p_1},\\
\end{align*}
where $r_2^{15} \left(H_{15}^{(1)}-p_1 \right), \ r_4^{15} \left(H_{15}^{(2)}-p_2 \right)$.

{\bf Hamiltonians $H_{150}^{(1)}=r_0^{15}\left(H_{15}^{(1)} \right), \ H_{150}^{(2)}=r_0^{15}\left(H_{15}^{(2)} \right)$}
\begin{align*}
&r_0^{150}:x_0=-(q_1 p_1-\beta_0)p_1, \ y_0=\frac{1}{p_1}, \ z_0=q_2, \ w_0=p_2, \\
&r_1^{150}:x_1=-(q_1 p_1-\beta_1-\beta_0)p_1, \ y_1=\frac{1}{p_1}, \ z_1=q_2, \ w_1=p_2, \\
&r_2^{150}:x_2=q_1+\frac{\beta_2-\beta_0}{p_1}+\frac{t-q_2}{p_1^2}, \ y_2=p_1, \ z_2=q_2, \ w_2=p_2-\frac{1}{p_1}, \\
&r_3^{150}:x_3=-(q_1 p_1-q_2 p_2-\beta_3-\beta_0-\beta_1-\beta_5)p_1, \ y_3=\frac{1}{p_1}, \ z_3=\frac{q_2}{p_1}, \ w_3=p_2 p_1, \\
&r_4^{150}:x_4=q_1, \ y_4=p_1, \ z_4=-((q_2-s)p_2-\beta_4)p_2, \ w_4=\frac{1}{p_2}, \\
&r_5^{150}:x_5=q_1+\frac{\beta_5 p_2}{p_1 p_2-1}, \ y_5=p_1, \ z_5=q_2+\frac{\beta_5 p_1}{p_1 p_2-1}, \ w_5=p_2,\\
&r_6^{150}:x_6=\frac{1}{q_1}, \ y_6=-(p_1 q_1-p_2 q_2+\beta_6-\beta_0)q_1, \ z_6=q_2 q_1, \ w_6=\frac{p_2}{q_1},\\
\end{align*}
where $r_2^{150} \left(H_{150}^{(1)}-\frac{1}{p_1} \right), \ r_4^{150} \left(H_{150}^{(2)}-p_2 \right)$.

{\bf Hamiltonians $H_{1503}^{(1)}=sr_3 \left(H_{150}^{(1)} \right), \ H_{1503}^{(2)}=sr_3 \left(H_{150}^{(2)} \right), \ sr_3:Q_1=q_1, \ P_1=p_1-\frac{q_2 p_2+\beta_3+\beta_0+\beta_1+\beta_5}{q_1}, \ Q_2=q_2 q_1, \ P_2=\frac{p_2}{q_1}$}
\begin{align}\label{reducedholo}
\begin{split}
&r_0^{1503}:x_0=-(q_1p_1+q_2p_2-\delta_0)p_1, \ y_0=\frac{1}{p_1}, \ z_0=q_2p_1, \ w_0=\frac{p_2}{p_1}, \\
&r_1^{1503}:x_1=-(q_1p_1+q_2p_2-\delta_1-\delta_0)p_1, \ y_1=\frac{1}{p_1}, \ z_1=q_2p_1, \ w_1=\frac{p_2}{p_1}, \\
&r_2^{1503}:x_2=q_1+\frac{q_2}{t}+\frac{2q_2p_2+2\delta_2}{p_1}+\frac{t}{p_1^2}, \ y_2=p_1, \ z_2=q_2p_1^2, \ w_2=\frac{p_2+\frac{p_1}{t}}{p_1^2}, \\
&r_3^{1503}:x_3=\frac{1}{q_1}, \ y_3=-(p_1q_1+\delta_3)q_1, \ z_3=q_2, \ w_3=p_2, \\
&r_4^{1503}:x_4=-\left( \left(q_1-\frac{1}{s} q_2 \right)p_1-\delta_4 \right)p_1, \ y_4=\frac{1}{p_1}, \ z_4=q_2, \ w_4=p_2+\frac{1}{s} p_1, \\
&r_5^{1503}:x_5=q_1p_2, \ y_5=\frac{p_1}{p_2}, \ z_5=q_2+\frac{q_1p_1+2\delta_5}{p_2}-\frac{1}{p_2^2}, \ w_5=p_2,\\
&r_6^{1503}:x_6=-(q_1p_1+q_2p_2-\delta_6-\delta_1-\delta_0)p_1, \ y_6=\frac{1}{p_1}, \ z_6=q_2 q_1, \ w_6=\frac{p_2}{q_1},
\end{split}
\end{align}
where $r_2^{1503} \left(H_{1503}^{(1)}+\frac{p_1 q_2}{t^2}-\frac{1}{p_1} \right), \ r_4^{1503} \left(H_{1503}^{(2)}+\frac{p_1 q_2}{s^2} \right)$.

Here, $\delta_j \ (j=0,1,\ldots,6)$ are complex constant parameters satisfying the parameter's relation\rm{: \rm}
\begin{align}
\begin{split}
&\beta_0=-\delta_6, \ \beta_1=-\delta_1, \ \beta_2=2\delta_0+2\delta_1+2\delta_2+\delta_6,\\
&\beta_3=-2\delta_0-\delta_1-2\delta_5-\delta_6, \ \beta_4=\delta_4, \ \beta_5=\delta_0+\delta_1+2\delta_5+\delta_6, \ \beta_6=\delta_0+\delta_1+\delta_3,
\end{split}
\end{align}
\begin{align}
\begin{split}
&3\delta_0+2\delta_1+2\delta_2+\delta_3+\delta_4+2\delta_5+\delta_6=1.
\end{split}
\end{align}
We remark that the transformations $r_2^{1503},r_5^{1503}$ are not its auto-B{\"a}cklund transformations. It is still an open question whether the transformations $r_2^{1503},r_5^{1503}$ can be considered as each transformation denoted by the symbol $\odot$ in the Oshima's paper (see \cite{Oshima}).

It is also still an open question whether we can obtain the Hamiltonian system with $H_{1503}^{(1)},H_{1503}^{(2)}$ by solving $3 \times 3$ Lax pair  (cf. \cite{KNS2,Oshima}) satisfying the following Riemann scheme:
\begin{equation*}
\begin{pmatrix}
X=0 & X=1 &  X=t &  X=\infty\\
\begin{matrix}
0 \\
\theta_1^{0} \\
\theta_2^{0}
\end{matrix}  & \begin{matrix}
0 \\
0 \\
\theta^{1}
\end{matrix}  &  \begin{matrix}
0 \\
0 \\
\theta^{t}
\end{matrix} &  \begin{matrix}
\theta_1^{\infty}\\
\theta_2^{\infty} \\
\theta_3^{\infty}
\end{matrix}
\end{pmatrix}
\end{equation*}

Here, we will conjecture the following relations between Riemann data and Holomorphy conditions $r_i^{1503} \ (i=0,1,\ldots,6);$

{\footnotesize
$\begin{pmatrix}
X=0 \\
 \begin{matrix}
0 \\
\delta_3 \\
\delta_4
\end{matrix}
\end{pmatrix} \Longleftrightarrow $ Holomorphy conditions $\begin{pmatrix}
\begin{matrix}
r_3^{1503} \\
r_4^{1503} 
\end{matrix}
\end{pmatrix}, \ \begin{pmatrix}
X=\infty \\
 \begin{matrix}
\delta_0 \\
\delta_1+\delta_0 \\
\delta_6+\delta_1+\delta_0
\end{matrix}
\end{pmatrix} \Longleftrightarrow $ Holomorphy conditions $\begin{pmatrix}
\begin{matrix}
r_0^{1503} \\
r_1^{1503} \\
r_6^{1503}
\end{matrix}
\end{pmatrix}.$

\vspace{0.5cm}

{\bf Appendix E: Holomorphy History}

{\footnotesize
\begin{tabular}{|c||c|c|} \hline
article & Author & Contents  \\ \hline
\cite{Pain1,Pain2} & P. Painlev\'e &  Convergence of meromorphic solution \\ \hline
\cite{O3} & K. Okamoto &  Patching data of space of initial conditions \\ \hline
\cite{Oka} & K. Okamoto and H. Kimura &  Patching data of Garnier system in n variables \\ \hline
\cite{T1,MMT} & A. Matumiya and K. Takano & Symplectic structure of space of initial conditions  \\ \hline
\cite{KK2,Suzuki} & H. Kimura and M. Suzuki & Degenerate Garnier System in two variables  \\ \hline
\cite{Ta} & N. Tahara &  Augmentation \\ \hline
\cite{YS} & Y. Yamada &  Relation between symmetry and holomorphy conditions \\ \hline
\end{tabular}
}

\end{document}